\documentclass[3p,12pt]{elsarticle}
\usepackage{color,xcolor}
\usepackage[colorlinks,linkcolor=blue,citecolor=blue]{hyperref}
\usepackage{amsmath,amsfonts,amsmath,amssymb,amsthm}
\usepackage{mathrsfs} 
\usepackage{graphicx,subcaption}
\usepackage{bm}
\usepackage[margin=10pt,font=small,labelfont=bf,labelsep=period]{caption}

\newcommand{\dx}{\,{\rm d}x}
\usepackage{multirow}
\usepackage{algorithm}
\usepackage{algpseudocode}

\numberwithin{equation}{section} 

\newtheorem{remark}{Remark}[section]

\biboptions{sort&compress}

\allowdisplaybreaks
\newtheoremstyle{mythm}{1.5ex plus 1ex minus .2ex}{1.5ex plus 1ex minus .2ex}
{\song}{\parindent}{\song\bfseries}{}{1em}{}
\theoremstyle{mythm}

\begin{document}

\begin{frontmatter}

\title{An efficient numerical method for surface acoustic wave equations over unbounded domains\footnote{The work was partially supported by the China National Key R\&D Project (Grant No. 2020YFA0709 800).}}
\date{Version 0.0, March 21, 2024}
\author{Jianguo Huang}
\ead{jghuang@sjtu.edu.cn}
\author{Likun Qiu\footnote{Corresponding author.}}
\ead{sjtumathqlk1997@sjtu.edu.cn}
\address{School of Mathematical Sciences, and MOE-LSC, Shanghai Jiao Tong University\\
Shanghai 200240, China}

\begin{abstract}
   Surface acoustic wave (SAW) devices are widely used in modern communication equipment and SAW equations describe the critical physical processes of acoustic-electric conversion in SAW devices. It is very challenging to numerically solve such equations, since they are typically three dimensional problems defined on unbounded domains. In this paper, we first use the perfectly matched layer method to truncate the unbounded domain and then propose a finite element tearing and interconnecting algorithm for the truncated equations based on the periodic structure of the truncated domain. We also design an effective solver for the ill-conditioned linear system of the Lagrange multipliers arising from discretization. Several numerical results are performed to demonstrate the efficiency of the proposed algorithm.
\end{abstract}


\end{frontmatter}
\section{Introduction}
SAW devices are important components of the radio-frequency (RF) and the intermediate-frequency (IF) stages in some communication systems, such as satellite receivers, remote control units, keyless entry systems, television sets, and mobile phones \cite{IntroSAW}. In the past 30 years, the rapid development of wireless communication has greatly stimulated the development of SAW devices. To achieve various functions, SAW devices have become increasingly complex. Therefore, efficient numerical simulation algorithms are crucial to the development of SAW devices. SAW equations are the essential mathematical and physical model that describe the acoustic-electric conversion in SAW devices and are included in almost all the mathematical models of SAW devices. In general, the finite element methods (FEMs) are used to solve SAW equations in some commercial simulation software due to their remarkable generality. However, solving SAW equations using the FEMs has been hampered by large memory consumption and low computational speed. To solve the SAW equations defined in large-scale domains, we need to design a more efficient algorithm with lower time and space complexity.

Generally, SAW devices are composed of piezoelectric substrate and interdigital electrodes. Commonly, we are only interested in the behavior of the regions beneath the electrodes, and the size of this region is much smaller than that of the whole device. Therefore, SAW equations are usually treated as semi-unbounded domain problems, with the SAW gradually decaying and vanishing along the piezoelectric substrate. The first step in solving the SAW equations is to address the semi-unbounded domain. As far as we know, there are several ways to tackle this problem. The spectral methods can solve this type of problem directly, since some of their basis functions are defined in unbounded domain and exhibit rapid decay (cf. \cite{Tissaoui-Kelly-Marras-2025,Shen-Tang-Wang-2011,Shen-Yu-2012,Shen-Yu-2010}). The boundary integral equation (BIE) methods are also commonly used to solve unbounded domain problems, convert a partial differential equation (PDE) defined in an unbounded domain into a boundary integral equation using Green's functions, and further develop numerical methods (cf. \cite{UnboundedBIE1,UnboundedBIE2,UnboundedBIE3,UnboundedBIE4}). The absorbing boundary conditions (ABC) and the perfectly matched layer (PML) methods are local methods, both of which truncate the unbounded domain by artificial boundary conditions. We refer the reader to \cite{ABC1,ABC2,ABC3,ABC4} for the ABC methods and \cite{PML1,ReviewOfPML,Leng-Ju-2019,Leng-Ju-2022} for the PML methods. Since the PML methods are simple in principle and easy to implement, they have been widely applied in many types of microacoustic devices, such as the SAW and the bulk acoustic wave (BAW) resonators (cf. \cite{PMLPiez1,PMLPiez2,PMLPiez3}). In this paper, we use the PML method to truncate the semi-unbounded domain.

So far, there have been several excellent numerical algorithms for SAW equations. In \cite{FEMandBEM2,FEMandBEM1}, the authors combined the finite element method and the boundary element method (BEM) to propose an efficient numerical method for solving the SAW equations with arbitrary geometries of the metallic electrodes. The periodic structure of the defined domain of the equations is a good starting point for designing efficient numerical algorithms. In \cite{HCT1,PMLPiez2}, the authors proposed the hierarchical cascading algorithm to calculate the admittance of SAW devices. It is a very efficient algorithm with low time and space complexity, but it performs poorly when we want to obtain all the numerical solutions within the region of interest.

The domain decomposition methods are highly suitable for solving equations defined in the domains with periodic structures \cite{DDM1}. Therefore, in this paper, we propose a finite element tearing and interconnecting (FETI) algorithm to solve the SAW equations truncated by the PML method. The FETI algorithm, first proposed by Farhat and Roux \cite{FETIOri}, is a type of the non-overlapping domain decomposition method. The introduction of the Dirichlet preconditioner and the lumped preconditioner \cite{preconditioner2FETI} greatly improves the efficiency of the FETI algorithm, making its convergence rate insensitive to the number of degrees of freedom (DOFs) in each subregion. In \cite{AnalysisOfFETI}, the author analyzed the convergence properties of the FETI algorithm in detail. Up to now, the FETI algorithm has been applied to the numerical solution of many partial differential equations. In \cite{FETIforDiffusionreaction}, the FETI algorithm was used to solve the diffusion-reaction equations. The Maxwell equations were solved by the FETI algorithm in \cite{FETIforMaxwell,FETIforDiscontinuousCoefficientsMaxwell}. In \cite{FETIBDNM} and \cite{FETIwithBdifferentiable}, the FETI algorithm was used to solve the 3D elastic frictional contact problems in combination with the B-differentiable Newton method (BDNM) and the B-differentiable equations (BDEs), respectively.

In the FETI algorithm, the efficiency of decoupling subregions (i.e., solving the Lagrange multiplier linear system) directly determines the overall efficiency of the algorithm. Generally, the Krylov subspace methods with the Dirichlet/lumped preconditioner are used in this step (cf. \cite{preconditioner2FETI}). Unfortunately, the linear system derived from the SAW equations is very ill-conditioned, and these  preconditioners do not perform well. Therefore, we propose a new method based on the Sherman-Morrison-Woodbury formula to solve the Lagrange multiplier linear system in this paper. The space complexity of this method is very low, and the computation time is within an acceptable range. Readers can find more details in Sect. \ref{Sect:Solve qusi-Toeplitz}, and we numerically demonstrate the efficiency of this method in Sect. \ref{sect:NumericalResult 3}.


We end this section by introducing some notations frequently used later on. For a bounded domain $\Omega\subset\mathbb{R}^d$ which has a Lipschitz boundary, let $L^2(\Omega)$ be the vector space consisting of all square integrable functions in $\Omega$, which is equipped with the norm
\[
\|f\|_{0,\Omega}:=\Big(\int_{\Omega}|f|^2 \dx\Big)^{\frac 1 2},\quad \forall f\in L^2(\Omega).
\]
Let $\alpha=(\alpha_1,...,\alpha_n)$ be a multi-index, where $\alpha_i$ is a non-negative integer, denote its length as $|\alpha|=\sum_{i=1}^{n}\alpha_i$, and write the weak partial derivative of a function $f$ with respect to $\alpha$ as
\[
D^{\alpha}f=\frac{\partial^{|\alpha|}}{\partial x_1^{\alpha_1}...\partial x^{\alpha_n}_n}f.
\]
Then, for any non-negative integer $s$, let $H^s(\Omega)$ be a Sobolev space consisting of all functions $f\in L^2(\Omega)$ such that $D^{\alpha}f\in L^2(\Omega)$ for all $|\alpha|\le s$.

The Frobenius norm $\|\cdot\|_F$ is defined as
$$\|\bm A\|_F:=\sqrt{{\rm tr}(\bm A\bm A^{\rm H})},\quad \forall \bm A\in\mathbb{C}^{n^2}.$$

 Throughout this paper, unless stated explicitly, we use Einstein's convention for summation whereby summation is implied when an index occurs exactly twice. For a vector with the $i$-th component $x_i$, we simply write $[x_i]$ or $\bm x$ as the vector itself. This convention also applies to any tensor.
\section{SAW equations}
\subsection{Mathematical model}
\begin{figure}[H]
	\centering
	\includegraphics[width=12cm]{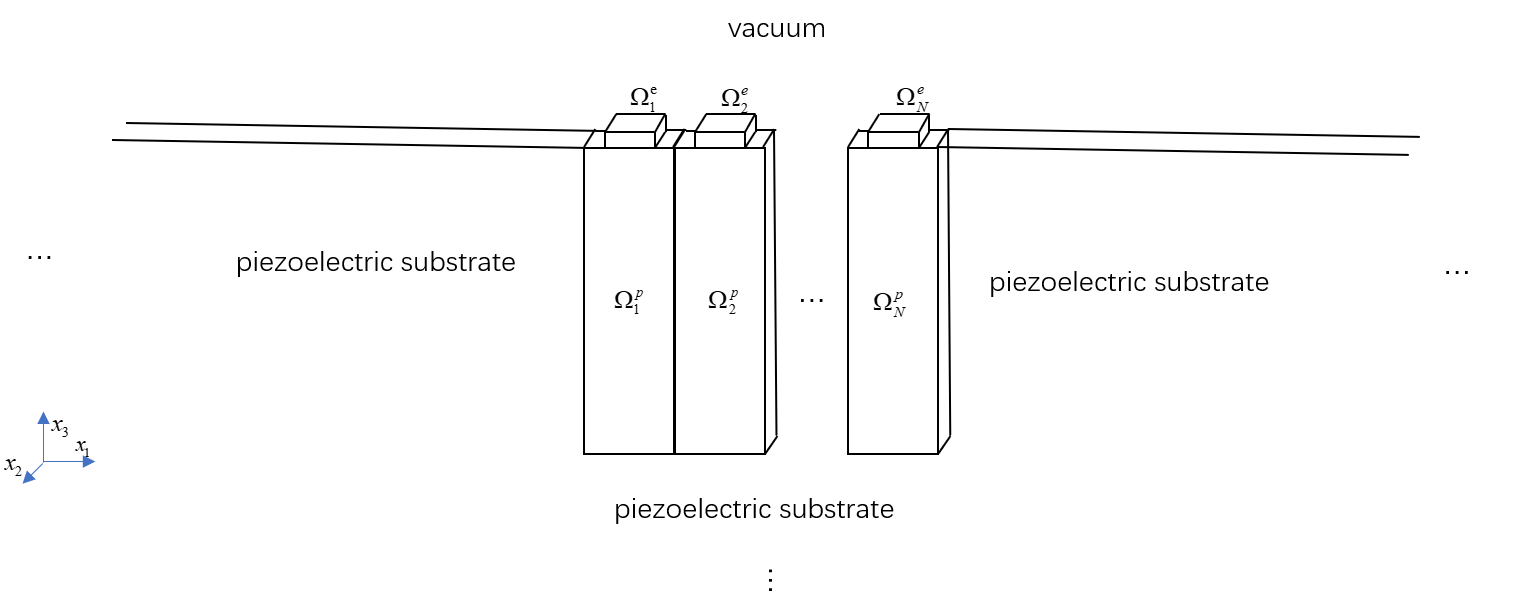}
	\caption{SAW device.}
	\label{fig:SAW device}
\end{figure}

Without loss of generality, we consider the SAW equations defined in the domain shown in Fig. \ref{fig:SAW device}, i.e. a  sequence of unit blocks
embedded in a large piezoelectric substrate, where the unit block is a pair of piezoelectric substrate $\Omega_i^p$ and electrode $\Omega_i^e$. When an alternating current electric input signal is applied to the electrodes, the electric field penetrates the piezoelectric substrate and SAW is induced due to the piezoelectric effect \cite{SAWMathModel,IntroSAW}.
 In general, we only want to obtain the numerical solutions in $\Omega=\cup_{i=1}^N(\Omega_i^{p}\cup \Omega_i^{e})$. The mathematical models can be described as the following partial differential equations.
\begin{equation}\label{eqn:SAW equations}
	\left\{
	\begin{array}{ll}
		\partial_j\sigma^p_{ij}+\omega^2\rho^pu_i=0&\text{ in }\Omega_m^{p},\\
		\partial_j\sigma^e_{ij}+\omega^2\rho^eu_i=0&\text{ in }\Omega_m^{e},\\
		\partial_jD_j^p=0&\text{ in }\Omega_m^{p}
	\end{array}\right.\quad (m=1,2,...,N),
\end{equation}
where $\bm \sigma^{\cdot}:=[\sigma_{ij}]$ is the stress tensor, ${\bm u}:=[u_j]$ is the displacement field, $\bm D^{\cdot}=[D_j]$ is the electric displacement vector, $\omega$ is the  frequency, $\rho^{\cdot}$ is the density, and the superscripts $p$ and $e$ represent the symbols correspond to the piezoelectric substrates and the electrodes, respectively. The constitutive relations in the piezoelectric substrates and the electrodes are:
$$
\begin{array}{ll}
	\left\{\begin{array}{l}
		\sigma_{ij}^p=\frac{1}{2}c_{ijkl}^p(\partial_ku_l+\partial_lu_k)+e_{kij}^p\partial_k\phi,\\
		D_i^p=\frac{1}{2}e_{ikl}^p(\partial_ku_l+\partial_lu_k)-\varepsilon_{ik}^p\partial_k\phi
	\end{array}\right.
	&\text{ in }\Omega_m^{p},\\
	\sigma^e_{ij}=\frac{1}{2}c_{ijkl}^e(\partial_ku_l+\partial_lu_k)&\text{ in }\Omega_m^{e},
\end{array}
$$
where $\phi$ is the electric potential, $[c_{ijkl}^{\cdot}]$ is the elasticity tensor, $[e_{kij}^{\cdot}]$ is the piezoelectric tensor and $[\varepsilon_{ik}^{\cdot}]$ is the dielectric permittivity tensor.

Let $\Gamma_{ij}^p:=\overline{\Omega_{i}^p}\cap\overline{\Omega_{j}^p}$ and $\Gamma_i:=\overline{\Omega_i^p}\cap \overline{\Omega_i^e}$. Each physical quantity satisfies the following continuity at the interfaces between the subregions.

$$
\left\{\begin{array}{lll}
	\bm \sigma^p|_{\Omega_i^p}\bm n_{i}^p+\bm \sigma^p|_{\Omega_j^p}\bm n_{j}^p=\bm 0,&\bm u|_{\Omega_i^p}=\bm u|_{\Omega_j^p},&\\
	\bm D|_{\Omega_i^p}\bm n_{i}^p+\bm D|_{\Omega_j^p}\bm n_{j}^p=0,&\phi|_{\Omega_i^p}=\phi|_{\Omega_j^p}&\text{ on }\Gamma_{ij}^p,\\
	\bm \sigma^p|_{\Omega_i^p}\bm n_{i}^p+\bm \sigma^e|_{\Omega_i^e}\bm n_{i}^e=\bm 0,& \bm u|_{\Omega_i^p}=\bm u|_{\Omega_i^e}&\text{ on }\Gamma_{i},
\end{array}\right.
$$
where $\bm n_i^p$ and $\bm n_i^e$ are the unit outward normal of the boundary face of $\Omega_i^p$ and $\Omega_i^e$.

We assume that the SAW attenuates to vanish along the positive and negative directions of the $x_1$-axis and the negative direction of the $x_3$-axis.  Alternating voltages $\phi_{0,i}$ are applied to each electrode $\Omega_{i}^e$. We can ignore the electric potential DOFs in all the electrodes, since the metal electrodes are equipotential bodies. In summary, the boundary conditions of the SAW equations are shown as follows.
$$
\left\{\begin{array}{ll}
	\lim_{x_1\to\pm\infty}\bm u=\lim_{x_3\to-\infty}\bm u=\bm 0,& \\
	\lim_{x_1\to\pm\infty}\phi=\lim_{x_3\to-\infty}\phi=0,& \\
	\phi=\phi_{0,i}&\text{ on }\Gamma_i,\\
	\bm \sigma^e\bm n=\bm \sigma^p\bm n=\bm 0,\quad \bm D^p\bm n=\bm D^e\bm n=0&\text{ on the rest boundary face,}
\end{array}\right.
$$
where $\bm n$ is the unit outward normal of the boundary face.
\subsection{Truncating the semi-unbounded domain by the perfectly matched layer method}
\begin{figure}[H]
	\centering
	\includegraphics[width=6cm]{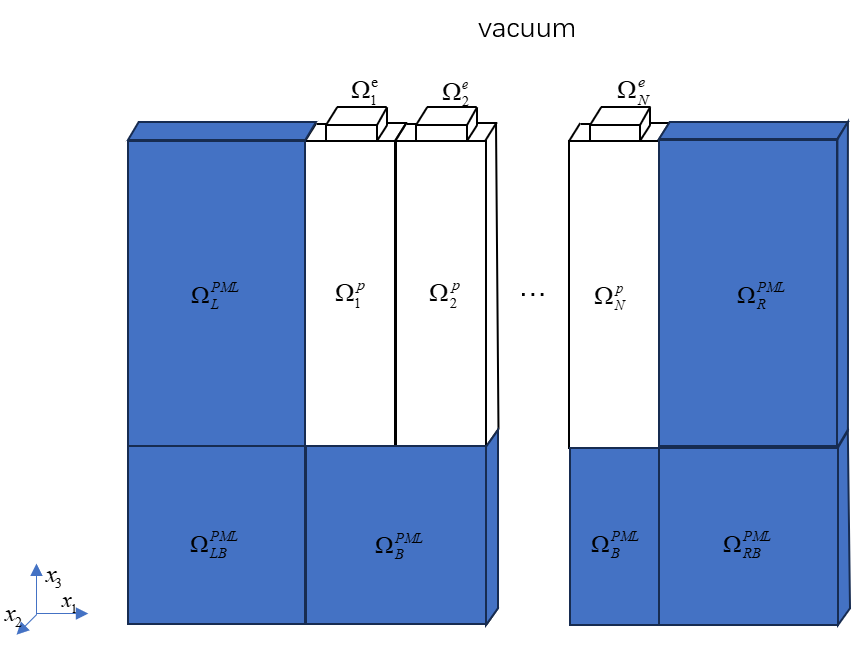}
	\caption{Truncted by PML.}
	\label{fig:SAW device PML}
\end{figure}

To solve the SAW equations numerically, we first need to truncate the semi-unbounded domain by the PML method, as shown in Fig. \ref{fig:SAW device PML}.
The SAW achieves exponential decay in the PML regions by the complex coordinate-stretching, which is defined as follows.
\begin{equation}\label{eqn:complex coordinate-stretching}
	\tilde{x}_k = x_k - {\rm i}\int^{x_k}_{0} d^{(m)}_k(\xi){\rm d\xi},\quad k=1,2,3,
\end{equation}
where ${\rm i}=\sqrt{-1}$ is the imaginary unit and $d^{(m)}_k(x_k)$ is the damping function of $\Omega^{\rm PML}_{m}$ ($m=$L,R,B,LB,LR), which is a continuous and monotonically
increasing positive-valued function with $x_k$ as variable (this property causes the solution to decay in the PML regions, see \cite{ReviewOfPML} for detail), such that

\begin{equation}\label{eqn:damping}
	\begin{cases}d^{(m)}_k\left(x_k\right)>0  \text{ in } \Omega^{\rm PML}_{m}, \\ d^{(m)}_k\left(x_k\right)=0  \text{ elsewhere. }\end{cases}
\end{equation}
For clarity, we denote the total PML region as $\Omega^{\rm PML}$, i.e. $\Omega^{\rm PML}:=\Omega^{\rm PML}_{\rm L}\cup\Omega^{\rm PML}_{\rm R}\cup\Omega^{\rm PML}_{\rm B}\cup\Omega^{\rm PML}_{\rm LB}\cup\Omega^{\rm PML}_{\rm RB}$
.  From (\ref{eqn:complex coordinate-stretching})-(\ref{eqn:damping}), we can see that the equations coincides with $(\ref{eqn:SAW equations})$ in $\Omega$ and the solution is continuous on $\Gamma_{\rm int}:=\overline{\Omega^{\rm PML}}\cap\overline{\Omega}$.  In ${\Omega}^{\rm PML}$ , the equations are shown as follows.
\begin{equation}\label{eqn:PML domain piez}
	\left\{\begin{array}{ll}
		\omega^2 \rho^p {u}_i+\frac{1}{\alpha_j}\partial_j \tilde{\sigma}_{ij}=0 &\text { in } {\Omega}^{\rm PML}, \\
		\frac{1}{\alpha_j}\partial_j\tilde{D}_j=0 &\text { in } {\Omega}^{\rm PML},\\
		{\bm u}=\bm 0\text{, }{\phi}=0 &\text{ on } \Gamma_{\rm ext},\\
		\tilde{\bm\sigma}\bm n=\bm 0\text{, }\tilde{\bm D}\bm n=0 &\text{ on } \partial{\Omega}^{\rm PML}\backslash(\Gamma_{\rm ext}\cup\partial\Omega),\\
			{\bm u}|_{\overline{\Omega^{\rm PML}}}=\bm u|_{\overline{\Omega}}\text{, }{\phi}|_{\overline{\Omega^{\rm PML}}}=\phi|_{\overline{\Omega}} &\text{ on }\Gamma_{\rm int},\\
			\tilde{\bm\sigma}\bm n^{\rm PML}+\bm{\sigma}^p\bm n^p=\bm 0, \tilde{\bm D}\bm n^{\rm PML}+\bm D^p\bm n^p=0&\text{ on }\Gamma_{\rm int},
		\end{array}\right.
	\end{equation}
where $\bm n^{\rm PML}$ and $\bm n^p$ are the unit outward normals of the boundary face of the PML region and the piezoelectric substrate region, $\Gamma_{\rm ext}$ is the set of the left, right and  bottom boundary faces of $\Omega^{\rm PML}$
and
$$
\alpha_k:=\frac{\partial\tilde{x}_k}{\partial x_k}=\left\{\begin{array}{ll}
	1-{\rm i} d_k^{(m)}(x_k) \text{ in } &\Omega_{m}^{\rm PML}\text{ $(m={\rm L, R, B, LB, RB})$ },\\
	1,&\text{ otherwise}.
\end{array}\right.
$$

The constitutive relations in PML region is
$$
\left\{\begin{array}{l}
	\tilde{\sigma}_{ij}=\frac{1}{2}c_{ijkl}^p(\frac{1}{\alpha_k}\partial_k{u}_l+\frac{1}{\alpha_l}\partial_l{u}_k)+e_{kij}^p\frac{1}{\alpha_k}\partial_k{\phi},\\
	\tilde{D}_j=\frac{1}{2}e_{ikl}^p(\frac{1}{\alpha_k}\partial_k{u}_l+\frac{1}{\alpha_l}\partial_l{u}_k)-\varepsilon_{ik}^p\frac{1}{\alpha_k}\partial_k{\phi}.
\end{array}\right.
$$

\subsection{Weak formulation for the SAW equations}
According to Fig. \ref{fig:SAW device PML}, the truncated domain is a combination of the following four types of subregion.
\begin{figure}[H]
	\centering
	\includegraphics[width=8cm]{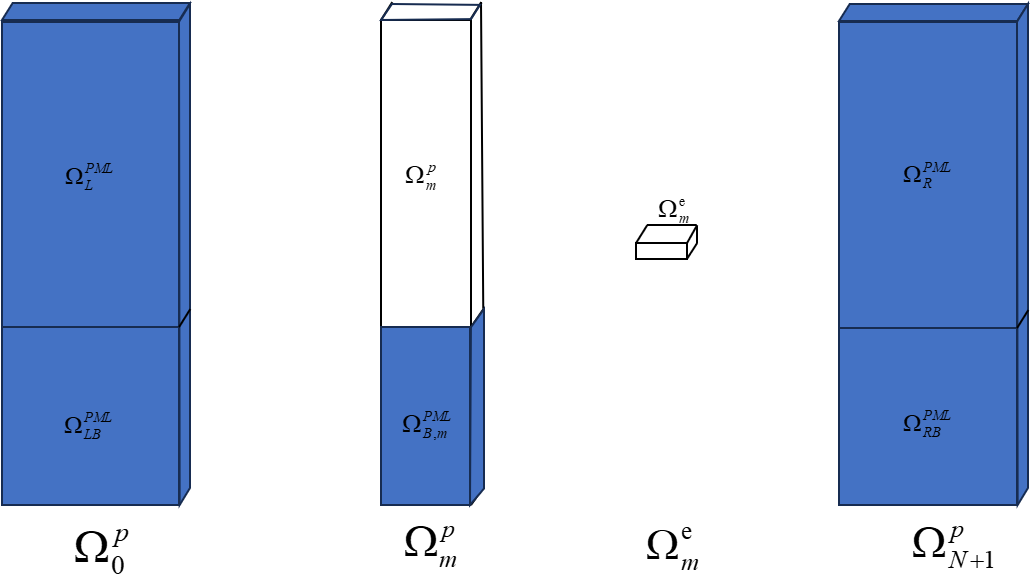}
	\caption{four types of subregion.}
	\label{fig:subregions}
\end{figure}
In Fig. \ref{fig:subregions}, $\Omega^{\rm PML}_{{\rm B},m}$ denotes the subregion of $\Omega_{\rm B}^{\rm PML}$ that is directly below $\Omega_{m}^p$. For convenience, we denote $\Omega^{\rm PML}_{\rm L}\cup \Omega^{\rm PML}_{\rm LB}$ as $\Omega_0^p$, $\Omega^{p}_m\cup \Omega^{\rm PML}_{{\rm B},m}$ as $\Omega^{p}_m$ and $\Omega^{\rm PML}_{\rm R}\cup \Omega^{\rm PML}_{\rm RB}$ as $\Omega_{N+1}^p$. In the remainder of this paper, we use the subscripts $l$, $t$, $r$ and $b$ to  denote the symbol associated with the left, top, right and bottom interfaces of a subregion, respectively.

Let
$$
\begin{aligned}
	&H_{\Gamma_{\rm ext}}^1(\Omega_m^p):=\{f:f\in H^1(\Omega_m^p)\text{ and } f=0\text{ on } \partial \Omega_m^p\cap\Gamma_{\rm ext} \},\\
	&H_{\Gamma_m}^1(\Omega^p_m):=\{f:f\in H^1(\Omega^p_m)\text{ and }f=0\text{ on } \Gamma_m\},\\
	&\bm u_m^p:=\bm u|_{\Omega_m^p},\quad \phi_m^p:=\phi|_{\Omega_m^p}(m=0,2,...,N+1),\\
	&\bm u_m^e:=\bm u|_{\Omega_m^e}(i=1,2,...,N),\\
\end{aligned}
$$
and $\bm\lambda_{m,i}^{\sigma}$, $\lambda_{m,i}^{D}$ $(i=l,t,r,b)$ be the Lagrange multipliers defined on the interface $i$ of $\Omega_m^p$, corresponding to the stress and electric displacement, respectively.

By the trace theorem, there are $\hat{\phi}_m^p\in L^{\frac{1}{2}}(\Omega_m^p)$ such that $\hat{\phi}_m^p|_{\Gamma_m}=\phi_{0,m}$ $(m=1,2,...,N)$. We denote $\bar{\phi}_m^p:={\phi}_m^p-\hat{\phi}_m^p$, and $\bar{\phi}|_{\Omega_m}:=\bar{\phi}_m^p$. Then, by integration by parts and transforming the integral region into the real domain, the weak formulation in each subregion can be expressed as follows.

\begin{itemize}
	\item \(\forall (\bm v,\psi)\in [H_{\Gamma_{ext}}^1(\Omega_{0}^p)]^3\times H_{\Gamma_{ext}}^1(\Omega_{0}^p)\):
	\begin{equation}\label{eqn:weakformulaPML_L}
		\left\{\begin{array}{l}
			A_0^p(\bm u_0^p,\bm v)+B_0^p(\bm v,\phi_0^p)=(\bm\lambda_{0,r}^{\sigma},\bm v)_{0,\Gamma_{0,1}^p},\\
			B_0^p(\bm u_0^p,\psi)-C_0^p(\phi_0^p,\psi)=(\lambda_{0,r}^{D},\psi)_{0,\Gamma_{0,1}^p}.
		\end{array}\right.
	\end{equation}
	\item \(\forall (\bm v,\psi)\in [H_{\Gamma_{ext}}^1(\Omega_{m}^p)]^3\times H_{\Gamma_m}^1(\Omega^p_m) \quad(m=1,2,...,N)\):
	\begin{equation}\label{eqn:weakformulaPiez}
		\left\{\begin{array}{ll}
			A_m^p(\bm u_m^p,\bm v)+B_m^p(\bm v,\bar{\phi}_m^p)=(\bm\lambda_{m,r}^{\sigma},\bm v)_{0,\Gamma_{m,m+1}^p}+(\bm\lambda_{m,t}^{\sigma},\bm v)_{0,\Gamma_{m}}-(\bm\lambda_{m-1,r}^{\sigma},\bm v)_{0,\Gamma_{m-1,m}^p}, \\
			B_m^p(\bm u_m^p,\psi)-C_m^p(\bar{\phi}_m^p,\psi)=(\lambda_{m,r}^{D},\psi)_{0,\Gamma_{m,m+1}^p}-(\lambda_{m-1,r}^{D},\psi)_{0,\Gamma_{m-1,m}^p}.
		\end{array}\right.
	\end{equation}
	\item \(\forall \bm v\in [H^1(\Omega_{m}^e)]^3 \quad (m=1,2,...,N)\):
	\begin{equation}
		A_m^e(\bm u_m^e,\bm v)=-(\bm\lambda_{m,t}^{\sigma},\bm v)_{0,\Gamma_m}.
	\end{equation}
	\item \(\forall (\bm v,\psi)\in [H_{\Gamma_{ext}}^1(\Omega_{N+1}^p)]^3\times H_{\Gamma_{ext}}^1(\Omega_{N+1}^p)\):
	\begin{equation}\label{eqn:weakformulaPML_R}
		\left\{\begin{array}{l}
			A_{N+1}^p(\bm u_{N+1}^p,\bm v)+B_{N+1}^p(\bm v,\phi_{N+1}^p)=-(\bm\lambda_{N,r}^{\sigma},\bm v)_{0,\Gamma_{N,N+1}^p},\\
			B_{N+1}^p(\bm u_{N+1}^p,\psi)-C_{N+1}^p(\phi_{N+1}^p,\psi)=-(\lambda_{N,r}^{D},\psi)_{0,\Gamma_{N,N+1}^p}.
		\end{array}\right.
	\end{equation}
\end{itemize}
and the continuity constraints on the interface:
\begin{equation}\label{eqn:weakformulaConti}
	\left\{\begin{aligned}
		&(\bm\mu,\bm u_m^p-\bm u_{m+1}^p)_{0,\Gamma_{m,m+1}^p}=\bm 0,\quad\forall \bm\mu\in [L^{\frac{1}{2}}(\Gamma_{m,m+1}^p)]^3\quad (m=0,...,N),\\
		&(\xi,\bar{\phi}_m^p-\bar{\phi}_{m+1}^p)_{0,\Gamma_{m,m+1}^p}=0,\quad\forall \xi\in L^{\frac{1}{2}}(\Gamma_{m,m+1}^p)\quad (m=0,...,N),\\
		&(\bm\mu,\bm u_m^p-\bm u_{m}^e)_{0,\Gamma_{m}}=\bm 0,\quad\forall \bm\mu\in [L^{\frac{1}{2}}(\Gamma_{m})]^3\quad (m=1,...,N),
	\end{aligned}\right.
\end{equation}
where $A_m^p(\bm u_m^p,\bm v)=A_{m,1}^p(\bm u_m^p,\bm v)-M_{m}^p(\bm u_m^p,\bm v)$ with
$$
\begin{aligned}
	&A_{m,1}^p(\bm u_m^p,\bm v)=\frac{1}{2}\int_{\Omega^{p}_m}c_{ijkl}^p\frac{\alpha_1\alpha_2\alpha_3}{\alpha_j}(\frac{1}{\alpha_k}\partial_k{u}_l+\frac{1}{\alpha_l}\partial_l{u}_k)\partial_j {v}_i {\rm dV} ,\\ &M_{m}^p(\bm u_m^p,\bm v)=\omega^2\rho^p\int_{\Omega^{p}_m}\alpha_1\alpha_2\alpha_3{u}_i{v}_i{\rm dV},
\end{aligned}
$$
$A_m^e(\bm u_m^e,\bm v)=A_{m,1}^e(\bm u_m^e,\bm v)-M_m^e(\bm u_m^e,\bm v)$ with
$$
\begin{aligned}
	&A_{m,1}^e(\bm u_m^e,\bm v)=\frac{1}{2}\int_{\Omega_m^e}c_{ijkl}^e(\partial_ku_l+\partial_lu_k)\partial_j v_i {\rm dV},\\
	&M_m^e(\bm u_m^e,\bm v)=\omega^2\rho^e\int_{\Omega_m^e} u_iv_i{\rm dV},
\end{aligned}
$$
and
$$
\begin{aligned}
	&B_m^p(\bm v,\bar{\phi}_m^p)=\int_{\Omega^{p}_m}e_{kij}^p\frac{\alpha_1\alpha_2\alpha_3}{\alpha_k\alpha_j}\partial_k\bar{\phi}\partial_j{v}_i{\rm dV},\quad C_m^p(\bar{\phi}_m^p,\psi)=\int_{\Omega^{p}_m} \varepsilon^p_{ik}\frac{\alpha_1\alpha_2\alpha_3}{\alpha_k\alpha_i}\partial_k\bar{\phi}\partial_i{\psi}{\rm dV},
\end{aligned}
$$
$(\cdot,\cdot)_{0,\Gamma}$ is the $L^2-$ inner product.
\section{Finite element tearing and interconnecting algorithm}\label{sect:FETI}

We discretize equations \eqref{eqn:weakformulaPML_L}-\eqref{eqn:weakformulaPML_R} by the Lagrange FEM and denote the matrices of  left hands as $\bm K_L$, $\bm K^p_m$, $\bm K^e_m$ and  $\bm K_R$ in sequence. We move the Dirichlet boundary conditions on $\Gamma_m$ to the right-hand and denote it as $\bm F^p_m$. We assume that the meshes of the same type of subregions are identical and that the grid points on the interfaces between any two subregions are matched. Therefore, the matrix of each $\Omega_{m}^p$ is the same, analogous to $\bm K^e_m$, we denote $\bm K_p:=\bm K^p_m$ and $\bm K_e:=\bm K^e_m$. Then, according to \cite{FETIOri}, the discrete forms of \eqref{eqn:weakformulaPML_L}-\eqref{eqn:weakformulaPML_R} are given as follows.
\begin{equation}\label{eqn:linearSystemSub-region}
	\left\{\begin{aligned}
		&\bm K_p\bm X_m^p=\bm F^p_m-\bm B_l^{p,\rm T}\bm\lambda_{m-1,r}+\bm B_t^{p,\rm T}\bm\lambda_{m,t}+\bm B_r^{p,\rm T}\bm\lambda_{m,r},&\\
		&\bm K_e\bm X_m^e=-\bm B_b^{e,\rm T}\bm\lambda_{m,t}&\quad (m=1,...N),\\
		&\bm K_L \bm X_L=\bm B_r^{L,\rm T}\bm\lambda_{0,r},\\
		&\bm K_R \bm X_{R}=-\bm B_l^{R,\rm T}\bm\lambda_{N,r},\\
	\end{aligned}\right.
\end{equation}
and the continuous conditions \eqref{eqn:weakformulaConti} on each interface are maintained by
\begin{equation}\label{eqn:continuous conditions}
	\left\{\begin{aligned}
		&\bm B_r^p \bm X^p_{m-1}=\bm B_l^p \bm X^p_{m}\quad (m=2,...,N-1),\\
		&\bm B_t^p \bm X^p_{m}=\bm B_b^e \bm X^e_{m}\quad (m=1,...,N), \\
		&\bm B_r^{L} \bm X_{L}=\bm B_l^p \bm X^p_{1},\\
		&\bm B_r^p \bm X^{p}_{N}=\bm B_l^{R} \bm X_{R}.
	\end{aligned}\right.
\end{equation}
In (\ref{eqn:linearSystemSub-region}) and (\ref{eqn:continuous conditions}), $\bm X^p_{m}$, $\bm X^e_{m}$, $\bm X_{L}$ and $\bm X_{R}$ are the numerical solutions, including the displacement and the electric potential, in $\Omega_{m}^p$,  $\Omega_{m}^e$ $(m=1,...,N)$, $\Omega_0^p$ and $\Omega_{N+1}^p$, respectively. $\bm\lambda_{m,i}$ is the discrete Lagrange multiplier vector defined on interface $i$ $(i=r, t)$ of the  subregion $\Omega_m^p$, containing both the Lagrange multipliers corresponding to the stress and the electric displacement. $\bm B_i^{j}$ is a $0$-$1$ matrix, which takes the following form.
$$
\bm B_i^{j}=\left[\begin{array}{ccccc}
	\bm 0&\cdots&\bm I&\cdots&\bm 0
\end{array}\right],
$$
where $\bm I$ is an identity matrix, whose size is equal to the number of DOFs on interface $i$. The number of columns of $\bm B_i^{j}$ is same as matrix $\bm K_j$.

Combining (\ref{eqn:linearSystemSub-region}) and (\ref{eqn:continuous conditions}), we can obtain the linear system for the Lagrange multipliers. We define the following matrices
\begin{equation}\label{eqn:matrixBlock}
	\begin{aligned}
	&\tilde{\bm A}_{rr}=\bm B^L_r\bm K_L^{-1}\bm B^{L,\rm T}_r,\\
	&\bm A_{ll} = \bm B^p_l \bm K_p^{-1}\bm B^{p,\rm T}_l,\quad \bm A_{lt} = -\bm B^p_l \bm K_p^{-1}\bm B^{p,\rm T}_t,\quad \bm A_{lr} = -\bm B^p_l \bm K_p^{-1}\bm B^{p,\rm T}_r,\\
	&\bm A_{tt} = \bm B^p_t \bm K_p^{-1}\bm B^{p,\rm T}_t+\bm B^e_b \bm K_e^{-1}\bm B^{e,\rm T}_b,\quad \bm A_{tr} = \bm B^p_t \bm K_p^{-1}\bm B^{p,\rm T}_r,\quad \bm A_{rr} = \bm B^p_r \bm K_p^{-1}\bm B^{p,\rm T}_r,\\
	&\tilde{\bm A}_{ll}=\bm B^R_l\bm K_R^{-1}\bm B^{R,\rm T}_l,
	\end{aligned}
\end{equation}
and the following vectors
\begin{equation}\label{eqn:vectorBlock}
	\bm b^p_{m,l} = \bm B^p_l\bm K_p^{-1}\bm F_m^p,\quad \bm b^p_{m,t} = -\bm B^p_t\bm K_p^{-1}\bm F_m^p,\quad \bm b^p_{m,r} = -\bm B^p_r\bm K_p^{-1}\bm F_m^p.
\end{equation}
Obviously, only applying different voltages to the electrodes $\Omega^e_i$ and $\Omega^e_j$ will cause $\bm F_i^p\neq \bm F_j^p$. Then, the linear system for the Lagrange multipliers is \begin{equation}\label{eqn:lagrangeMultiplierSystem}
	\bm A\bm\lambda=\bm b,
\end{equation} where
\begin{equation}
	\begin{aligned}
		&\bm A=\left[\begin{array}{ccccccccc}
			\tilde{\bm A}_{rr}+\bm A_{ll} & \bm A_{lt} & \bm A_{lr} & & & & & \\
			\bm A_{lt}^{\rm T} & \bm A_{tt} & \bm A_{tr} & & & & & \\
			\bm A_{lr}^{\rm T} & \bm A_{tr}^{\rm T} & \bm A_{rr}+\bm A_{ll} & \bm A_{lt} & \bm A_{lr} & & & \\
			& & \bm A_{lt}^{\rm T} & \bm A_{tt} & \bm A_{tr} & & & \\
			& & \bm A_{lr}^{\rm T} & \bm A_{tr}^{\rm T} & \bm A_{rr}+\bm A_{ll} & & & \\
			& & & & & \ddots & & & \\
			& & & & & & \bm A_{rr}+\bm A_{ll} & \bm A_{lt} & \bm A_{lr} \\
			& & & & & & \bm A_{lt}^{\rm T} & \bm A_{tt} & \bm A_{tr} \\
			& & & & & & \bm A_{lr}^{\rm T} & \bm A_{tr}^{\rm T} & \bm A_{rr}+ \tilde{\bm A}_{ll}
		\end{array}\right],\\
	&\bm\lambda = \left[\begin{array}{c}
		\bm\lambda_{0,r}\\
		\bm\lambda_{1,t}\\
		\bm\lambda_{1,r}\\
		\bm\lambda_{2,t}\\
		\bm\lambda_{2,r}\\
		...\\
		\bm\lambda_{N-1,r}\\
		\bm\lambda_{N,t}\\
		\bm\lambda_{N,r}
	\end{array}\right],\quad
	\bm b=\left[\begin{array}{c}
		\bm b^p_{1,l}\\
		\bm b^p_{1,t}\\
		\bm b^p_{1,r}+\bm b^p_{2,l}\\
		...\\
		\bm b^p_{N-1,r}+\bm b^p_{N,l}\\
		\bm b^p_{N,t}\\
		\bm b^p_{N,r}
	\end{array}\right].
	\end{aligned}
\end{equation}
Once we get $\bm\lambda$, all the subregions are decoupled. We can obtain the numerical solutions in each subregion through the following process.
\begin{equation}\label{eqn:numerical solution}
	\begin{aligned}
		&\bm X_m^p=\bm K_p^{-1}\bm F_m^p-\bm K_p^{-1}\bm B_l^{p,\rm T}\bm\lambda_{m-1,r}+\bm K_p^{-1}\bm B_t^{p,\rm T}\bm\lambda_{m,t}+\bm K_p^{-1}\bm B_r^{p,\rm T}\bm\lambda_{m,r},\\
		&\bm X_m^e=-\bm K_e^{-1}\bm B_b^{e,\rm T}\bm\lambda_{m,t}\quad (m=1,2,...,N).
	\end{aligned}
\end{equation}
\begin{remark}
	The matrix $\bm K_e$ seems to be singular since the boundary conditions are not sufficient to guarantee an unique solution in $\Omega_m^e$. However, in fact, it is nonsingular when $\omega>0$. We denote the matrices obtained by discretizing $A_m^e(\bm u_m^e,\bm v)$ and $M_m^e(\bm u_m^e,\bm v)$ as $\bm K_{h}$ and $\bm M_{h}$, respectively. Then, we have $\bm K_e=\bm K_{h}-\bm M_{h}$. Clearly, $\bm K_{h}$ is singular, $\bm M_{h}$ is nonsingular, and therefore $\bm K_e$ is nonsingular.
\end{remark}
\begin{remark}
	The linear system \eqref{eqn:lagrangeMultiplierSystem} can be easily solved in parallel using Krylov subspace methods if $A$ is positive definite or not too ill-conditioned. Unfortunately, the linear system derived from the SAW equations is indefinite and highly ill-conditioned, making the classical Dirichlet preconditioner and lumped preconditioner (cf. \cite{preconditioner2FETI}) inefficient. Therefore, we propose an efficient direct method in the next section.
\end{remark}
\begin{remark}
	In \eqref{eqn:matrixBlock}, \eqref{eqn:vectorBlock} and \eqref{eqn:numerical solution}, we need to solve $6+m_{\phi}$  linear systems ($m_{\phi}$ is the number of  different $\phi_{0,m}$), whose sizes only depend on the mesh of the subregions. These linear systems are $\bm K_L^{-1}\bm B_{r}^{L,\rm T}$, $\bm K_p^{-1}\bm B_{l}^{p,\rm T}$, $\bm K_p^{-1}\bm B_{t}^{p,\rm T}$, $\bm K_p^{-1}\bm B_{r}^{p,\rm T}$, $\bm K_e^{-1}\bm B_{b}^{e,\rm T}$, $\bm K_R^{-1}\bm B_{l}^{R,\rm T}$ and $\bm K_p^{-1}\bm F_m^p$. All of these can be computed in parallel, and are independent of $N$. When the mesh of the subregions is not changed but $N$ is increased, we can store the solutions of these linear systems and reuse them.
	
\end{remark}
\section{Efficiently and accurately implement the FETI algorithm}
In the previous section, we introduced the process of the FETI algorithm for solving the SAW equations. Since the linear systems derived from the SAW equations are typically highly ill-conditioned, special treatments are required to ensure the efficiency and accuracy of the algorithm. In the following two subsections, we introduce these two special treatments separately.
\subsection{Dimensionless}
In practical applications, the orders of magnitude of the material parameters vary significantly, which can cause algorithm instability. For example, the orders of the magnitude of $c_{ijkl}^p$ and $\varepsilon_{ik}^p$ are approximately $10^{10}$ and $10^{-12}$, respectively, when we solve the SAW equations with the YX128LN piezoelectric substrate. Therefore, we need to do some preprocessing before the algorithm starts. For convenience, we denote the matrices discretized from $A_{m,1}^p(\bm u_m^p,{\bm v})$, $M_m^p(\bm u_m^p,\bm v)$, $A_{m,1}^e(\bm u_m^p,{\bm v})$, $M_m^e(\bm u_m^p,\bm v)$, $B_m^p(\bm v,\bar{\phi}_m^p)$, $C_m^p(\bar{\phi}_m^p,\psi)$ as $\bm K_{uu}^p$, $\bm M_{uu}^p$, $\bm K_{uu}^e$, $\bm M_{uu}^e$,  $\bm K_{up}$ and $\bm K_{pp}$, respectively. The goal of the dimensionless process is to scale the parameters so that the six matrices have a similar order of magnitude.  We require the dimensionless parameters $\bar{c}_{ijkl}^p$, $\bar{c}_{ijkl}^e$, $\bar{\omega}$, $\bar{e}_{kij}^p$, $\bar{\varepsilon}_{ik}^p$, $\bar{\rho}^p$, $\bar{\rho}^e$ and $\bar{\bm x}$ satisfy
\begin{equation}\label{eqn:Dimensionless}
\begin{aligned}
		&{c}_{ijkl}^p=c_1\bar{c}_{ijkl}^p,\quad {c}_{ijkl}^e=c_1\bar{c}_{ijkl}^e,\quad \omega=\omega_1\bar{\omega},\quad e_{kij}^p=e_1\bar{e}_{kij}^p,\quad \varepsilon_{ik}^p=\varepsilon_1\bar{\varepsilon}_{ik}^p,\\ &\rho^p=\rho_1\bar{\rho}^p,\quad \rho^e=\rho_1\bar{\rho}^e,\quad \bm x=l_1\bar{\bm x},
\end{aligned}
\end{equation}
where $c_1$, $\omega_1$, $e_1$, $\varepsilon_1$, $\rho_1$ and $l_1$ are some constants satisfy
\begin{equation}\label{eqn:PropOftheDimensionlessConstant}
	c_1=\varepsilon_1^{-1},\quad l_1=\sqrt{\frac{c_1}{\omega_1^2\rho_1}}\text{ and }\quad e_1=1. 
\end{equation}
Then, we replace the original parameters by the dimensionless ones in \eqref{eqn:weakformulaPML_L}-\eqref{eqn:weakformulaPML_R} and denote the matrices discretized from the ``dimensionless equations" as $\bar{\bm K}_{uu}^p$, $\bar{\bm M}_{uu}^p$, $\bar{\bm K}_{uu}^e$, $\bar{\bm M}_{uu}^e$, $\bar{\bm K}_{up}$ and $\bar{\bm K}_{pp}$, which satisfy
\begin{align}
	&\left[\begin{array}{cc}
		\bm K_{uu}^p-\bm M_{uu}^p & \bm K_{up}\\
		\bm K_{up}^{\rm T} & -\bm K_{pp}
	\end{array}\right]=\left[\begin{array}{cc}
	\frac{c_1}{l_1^2}\bar{\bm K}_{uu}^p-\omega_1^2\rho_1\bar{\bm M}_{uu}^p & \frac{e_1}{l_1^2}\bar{\bm K}_{up}\\
	\frac{e_1}{l_1^2}\bar{\bm K}_{up}^{\rm T} & -\frac{\varepsilon_1}{l_1^2}\bar{\bm K}_{pp}
\end{array}\right],\\
& \bm K_{uu}^e-\bm M_{uu}^e=\frac{c_1}{l_1^2}\bar{\bm K}_{uu}^e-\omega_1^2\rho_1\bar{\bm M}_{uu}^e.
\end{align}
Since (\ref{eqn:PropOftheDimensionlessConstant}) holds, the linear systems
$$
\left[\begin{array}{cc}
\bm K_{uu}^p-\bm M_{uu}^p & \bm K_{up}\\
\bm K_{up}^{\rm T} & -\bm K_{pp}
\end{array}\right]\left[\begin{array}{l}
\bm u\\ \phi
\end{array}\right]=\bm F
$$
and
$$
(\bm K_{uu}^e-\bm M_{uu}^e)\bm u=\bm 0
$$
are equivalent to
\begin{equation}\label{eqn:pre-system}
	\left[\begin{array}{cc}
	\bar{\bm K}_{uu}^p-\bar{\bm M}_{uu}^p & \bar{\bm K}_{up}\\
	\bar{\bm K}_{up}^{\rm T} & -\bar{\bm K}_{pp}
	\end{array}\right]
	\left[\begin{array}{l}
	\bar{\bm u}\\ \bar{\phi}
	\end{array}\right]=\left[\begin{array}{cc}
	a\bm I & \bm 0\\
	\bm 0 & b\bm I
	\end{array}\right]^{-1}\bm F,
\end{equation}
and
\begin{equation}\label{eqn:pre-system2}
(\bar{\bm K}_{uu}^e-\bar{\bm M}_{uu}^e)\bar{\bm u}=\bm 0,
\end{equation}
where
$a=\sqrt{c_1}/l_1$, $b=\sqrt{{\varepsilon_1}}/{l_1}$, $\bar{\bm u}=a\bm u$ and $\bar{\phi}=b\phi$. Therefore, it is very convenient to restore the solutions of (\ref{eqn:pre-system})-\eqref{eqn:pre-system2} to the original solution.

Through the above process, the orders of the magnitude of the matrix blocks $\bar{\bm K}_{uu}^p$, $\bar{\bm M}_{uu}^p$, $\bar{\bm K}_{uu}^e$, $\bar{\bm M}_{uu}^e$, $\bar{\bm K}_{up}$ and $\bar{\bm K}_{pp}$ become similar. Therefore, the linear system derived from \eqref{eqn:weakformulaPML_L}-\eqref{eqn:weakformulaPML_R} becomes more stable. Since the FETI algorithm is based on \eqref{eqn:pre-system} and \eqref{eqn:pre-system2}, it will also be more stable.
\subsection{Efficiently solving the linear system for the Lagrange multipliers }\label{Sect:Solve qusi-Toeplitz}
The most important part of the FETI algorithm is solving (\ref{eqn:lagrangeMultiplierSystem}). The size of this linear system is much smaller than that of the discrete SAW equations for a single unit block when $N$ is small. However, when $N$ becomes very large, solving it requires a significant amount of storage and time. Therefore, it is crucial to find an efficient algorithm. Unfortunately, the linear system is indefinite and highly ill-condition, making the classical Dirichlet preconditioner and lumped preconditioner (c.f. \cite{preconditioner2FETI}) inefficient. In this section, we propose an algorithm with low space complexity and time complexity, based on the periodic structure of the domain shown in Fig. \ref{fig:SAW device PML}.

According to the form of the matrix $\bm A$ in (\ref{eqn:lagrangeMultiplierSystem}), we can convert $\bm A$ into a block tridiagonal qusi-Toeplitz matrix by the following process.

Let
\begin{equation}\label{eqn:MatrixBlocks}
	{\begin{aligned}
	&\bm M_L=\left[\begin{array}{cc}
	s\bm I&\bm 0\\
	\bm 0&\tilde{\bm A}_{rr}+\bm A_{ll}
\end{array}\right]
,\quad
\bm M=\left[\begin{array}{cc}
	\bm A_{tt}&\bm A_{tr}\\
	\bm A_{tr}^{\rm T}&\bm A_{rr}+\bm A_{ll}
\end{array}\right]
,\quad
\bm M_R=\left[\begin{array}{cc}
	\bm A_{tt}&\bm A_{tr}\\
	\bm A_{tr}^{\rm T}&\bm A_{rr}+\tilde{\bm A}_{ll}
\end{array}\right],\\
&\bm B=\left[\begin{array}{cc}
	\bm 0&\bm A_{lt}^{\rm T}\\
	\bm 0&\bm A_{lr}^{\rm T}
\end{array}\right],\quad
\bm \lambda_0=\left[\begin{array}{cc}
	\tilde{\bm \lambda}\\
	\bm \lambda_{0,r}
\end{array}\right],\quad
\bm\lambda_i=\left[\begin{array}{cc}
	\bm\lambda_{i,t}\\
	\bm\lambda_{i,r}
\end{array}\right],\quad
\bm b_L=\left[\begin{array}{cc}
	\bm 0\\
	\bm b_{1,l}^p
\end{array}\right],\quad
\bm b_i=\left[\begin{array}{cc}
	\bm b_{i,t}\\
	\bm b_{i,r}+\bm b_{i,l}
\end{array}\right],
\end{aligned}}
\end{equation}
where $\tilde{\bm \lambda}$ is an auxiliary variable and $s=\|\bm \tilde{\bm A}_{rr}+\bm A_{ll}\|_{F}$. The constant $s$ can prevent $\bm M_L$ from being singular. Then, the linear system (\ref{eqn:lagrangeMultiplierSystem}) is equivalent to the following system.
\begin{equation}\label{eqn:lagrangeSystemBlock}
	\left[\begin{array}{ccccc}
		\bm M_L & \bm B^{\mathrm{T}} & & & \\
		\bm B & \bm M & \bm B^{\mathrm{T}} & & \\
		& \bm B & \ddots & \ddots & \\
		& & \ddots & \bm M & \bm B^{\mathrm{T}} \\
		& & & \bm B & \bm M_R
	\end{array}\right]
	\left[\begin{array}{c}
		\bm\lambda_0\\
		\bm\lambda_1\\
		...\\
		\bm\lambda_{N-1}\\
		\bm\lambda_{N}
	\end{array}\right]
	=\left[\begin{array}{c}
		\bm b_L\\
		\bm b_1\\
		...\\
		\bm b_{N-1}\\
		\bm b_{N,t}\\
		\bm b_{N,r}
	\end{array}\right].
\end{equation}
For convenience, we denote (\ref{eqn:lagrangeSystemBlock}) as $\tilde{\bm A}\tilde{\bm\lambda}=\tilde{\bm F}$. $\tilde{\bm A}$ is a block tridiagonal quasi-Toeplitz matrix, i.e., a Toeplitz matrix with low rank perturbation. There are two main ideas to solve this type of linear system, the first is based on the  Sherman-Morrison-Woodbury formula (c.f.\cite{Qusi-Toeplitz2,Qusi-Toeplitz3}) and the other is called cyclic reduction (c.f.\cite{Qusi-Toeplitz1}). In this paper, we choose the former way. We define matrices $\bm L$ and $\bm \Lambda$ as
$$
\bm L=\left[\begin{array}{cccc}
	\bm I& & & \\
	\bm L_1& \bm I& & \\
	 &\ddots& \ddots& \\
	 & & \bm L_1&\bm I
\end{array}\right]\text{ and }
\bm\Lambda=\left[\begin{array}{cccc}
	\bm\Lambda_1& & & \\
	& \ddots& & \\
	& &\bm\Lambda_1 & \\
	& & &\bm\Lambda_2
\end{array}\right].
$$
Then,
$$
\bm L\bm\Lambda\bm L^{\rm T}=\left[\begin{array}{ccccc}
	\bm\Lambda_1&\bm\Lambda_1\bm L_1^{\rm T} & & & \\
	\bm L_1\bm\Lambda_1&\bm L_1\bm\Lambda_1\bm L_1^{\rm T}+\bm\Lambda_1&\bm\Lambda_1\bm L_1^{\rm T} & & \\
	&\ddots&\ddots&\ddots & \\
	& &\bm L_1\bm\Lambda_1&\bm L_1\bm\Lambda_1\bm L_1^{\rm T}+\bm\Lambda_1&\bm\Lambda_1\bm L_1^{\rm T}\\
	& & &\bm L_1\bm\Lambda_1&\bm L_1\bm\Lambda_1\bm L_1^{\rm T}+\bm\Lambda_2
\end{array}\right].
$$
We require the matrices $\bm L_1$, $\bm \Lambda_1$ and $\bm \Lambda_2$ satisfy
\begin{equation}\label{eqn:matrixEquations}
	\left\{\begin{array}{l}
		\bm L_1\bm\Lambda_1=\bm B,\\
		\bm L_1\bm\Lambda_1\bm L_1^{\rm T}+\bm\Lambda_1=\bm M,\\
		\bm L_1\bm\Lambda_1\bm L_1^{\rm T}+\bm\Lambda_2=\bm M_R.
	\end{array}\right.
\end{equation}
The first two equations in (\ref{eqn:matrixEquations}) are equivalent to
\begin{equation}\label{eqn:main nonlinear system}
	\bm B\bm\Lambda_1^{-1}\bm B^{\rm T}+\bm\Lambda_1=\bm M.
\end{equation}
If the matrices $\bm B$ and $\bm M$ are well-condition, (\ref{eqn:main nonlinear system}) can be solved rapidly by the Newton's method (c.f.\cite{Newton'sMethodForNonlinear}). However, both of $\bm B$ and $\bm M$ deduced from the SAW equations are sick. Let's introduce the best treatment we've tested, which is a combination of the Double method \cite{DoubleMethod1} and the Newton's method for quadratic matrix equations \cite{NewtonMethodForQuadratric}. Obviously, the solutions of (\ref{eqn:main nonlinear system}) can solve the following quadratic matrix equation:
\begin{equation}\label{eqn:quadratic matrix equation}
	Q(\bm Y):=-\bm B^{\rm T}+\bm M\bm Y-\bm B\bm Y^2=\bm 0,
\end{equation}
where $\bm Y=\bm\Lambda_1^{-1}\bm B^{\rm T}$. If we use the Newton's method to solve (\ref{eqn:quadratic matrix equation}) directly with an arbitrary initial solution, it may converge to  (\ref{eqn:quadratic matrix equation})'s own solution, which may not be able to solve (\ref{eqn:main nonlinear system}). Therefore, we compute an nice initial solution by \eqref{eqn:main nonlinear system}'s Double method first,  and then solve it by \eqref{eqn:quadratic matrix equation}'s Newton's method.

The relative residuals $\rho_{D}$ and $\rho_{N}$ are defined as follows.
$$
\begin{aligned}
	&\rho_D(\bm\Lambda_1^{k+1},\bm\Lambda_1^{k})=\frac{\|\bm\Lambda_1^{k+1}-\bm\Lambda_1^{k}\|_F}{\|\bm\Lambda_1^{k}\|_F},\\
	&\rho_{N}(\bm Y^k)=\frac{\|Q(\bm Y^k)\|_F}{\|\bm B\|_F\|\bm Y^k\|_F^2+\|\bm M\|_F\|\bm Y^k\|_F+\|\bm B^{\rm T}\|_F}.
\end{aligned}
$$
Then, the algorithm can be described as follows.
\begin{algorithm}[H]
	\caption{Double-Newton method}\label{alg:Double-Newton}
	\begin{algorithmic}[1]
		\State input: $\epsilon_D$, $\epsilon_N$ and $Iter_{max}$
		\State set $\bm B^0=\bm B^{\rm T}$, $\bm \Lambda_1^0=\bm M$ and $\bm P^0=\bm 0$\Comment{The Double method begins}
		\For{$k=1,2,3...$}
		\State $\bm B^{k}=\bm B^{k-1}(\bm \Lambda_1^{k-1}-\bm P^{k-1})^{-1}\bm B^{k-1}$
		\State $\bm \Lambda_1^{k}=\bm \Lambda_1^{k-1}-\bm B^{k-1, \rm T}(\bm \Lambda_1^{k-1}-\bm P^{k-1})^{-1}\bm B^{k-1}$
		\State $\bm P^{k}=\bm P^{k-1}+\bm B^{k-1}(\bm \Lambda_1^{k-1}-\bm P^{k-1})^{-1}\bm B^{k-1,\rm T}$
		\If{$\rho_D(\bm\Lambda_1^{k},\bm\Lambda_1^{k-1})<\epsilon_D$}
		\State $\bm\Lambda_{1*}=\bm\Lambda_{1}^k$
		\State break.
		\EndIf
		\EndFor
		\State set $\bm Y^0=\bm\Lambda_{1*}^{-1}\bm B^{\rm T}$, $k=0$\Comment{The Newton's method begins}
		\While{$\rho_{N}(\bm Y^k)>\epsilon_N$ and $k<Iter_{max}$}
		\State solve $\bm B \bm E_N\bm Y^k+(\bm B\bm Y^k-\bm M)\bm E_N=Q(\bm Y^k)$
		\State $\bm Y^{k+1}=\bm Y^{k} +\bm E_N$
		\State $k=k+1$
		\EndWhile
		\State $\bm\Lambda_{1*}=\bm M-\bm B\bm Y^k$
		\State output: $\bm\Lambda_{1*}$
	\end{algorithmic}
\end{algorithm}
\begin{remark}
	The generalized Sylvester equation in line 14 of  Algorithm \ref{alg:Double-Newton} can be solved directly using the Schur decomposition and the Hessenberg-triangular decomposition (c.f.\cite{MatrixComputations}). Let $n_m$ be the size of $\bm M$, the time complexity of Algorithm \ref{alg:Double-Newton} is $O((n_D+n_N)n_m^3)$, where $n_D$ and $n_N$ are the number of iterations of the double method and Newton's method, respectively. Obviously, the time complexity is not related to $N$.
\end{remark}
Now, we obtain a decomposition of $\tilde{\bm A}$:
\begin{equation}
	\tilde{\bm A}=\bm L\bm \Lambda\bm L^{\rm T}+\bm E_1\bm M_1^{\rm T},
\end{equation}
where $\bm E_1=\left[\begin{array}{llll}
	\bm I & \bm 0 & \cdots & \bm 0
\end{array}\right]^{\rm T}$ and $\bm M_1^{\rm T}=\left[\begin{array}{llll}
\bm M_L-\bm \Lambda_{1} & \bm 0 & \cdots & \bm 0
\end{array}\right]$. Thus, the solution $\tilde{\bm\lambda}$ is obtained by the Sherman-Morrison-Woodbury formula:
\begin{equation}\label{eqn:Sherman-Morrison-Woodbury formula}
	\begin{aligned}
		\tilde{\bm\lambda}&=(\bm L\bm \Lambda\bm L^{\rm T}+\bm E_1\bm M_1^{\rm T})^{-1}\tilde{\bm F}\\
		&=(\bm L\bm \Lambda\bm L^{\rm T})^{-1}\tilde{\bm F}-(\bm L\bm \Lambda\bm L^{\rm T})^{-1}\bm E_1(\bm I+\bm M_1^{\rm T}(\bm L\bm \Lambda\bm L^{\rm T})^{-1}\bm E_1)^{-1}\bm M_1^{\rm T}(\bm L\bm \Lambda\bm L^{\rm T})^{-1}\tilde{\bm F}.
	\end{aligned}
\end{equation}
We can compute (\ref{eqn:Sherman-Morrison-Woodbury formula}) efficiently by the following algorithm.
\begin{algorithm}[H]\caption{compute (\ref{eqn:Sherman-Morrison-Woodbury formula})}\label{alg:Sherman-Morrison-Woodbury formula}
	\begin{algorithmic}[1]
			\State let $\bm F:=[\tilde{\bm F},\bm E_1]$
		\For{$i=2,...,N+1$}
		\State $\bm F_i:=\bm F_i-\bm L_1\bm F_{i-1}$\Comment{{\scriptsize $\bm F_i$ represents the rows $(i-1)n_m+1$ to $in_m$ of $\bm F$.} }
		\EndFor
		\State let $\hat{\bm F}_1:=[\bm F_2,\bm F_3,...,\bm F_N]$
		\State $\hat{\bm F}_1:=\bm\Lambda_{1}^{-1}\hat{\bm F}_1$
		\State $\bm F_{N+1}:=\bm \Lambda_2^{-1}\bm F_{N+1}$
		\State let $\bm F:=[\hat{\bm F}_{1,1}^{\rm T},\hat{\bm F}_{1,2}^{\rm T},...,\hat{\bm F}_{1,N}^{\rm T},\bm F_{N+1}^{\rm T}]^{\rm T}$\Comment{{\scriptsize $\hat{\bm F}_{1,i}$ represents the columns $(i-1)n_m+1$ to $in_m$ of $\hat{\bm F}_1$.}}
		\For{$i=N,N-1,...,1$}
		\State $\bm F_i:=\bm F_i-\bm L_1^{\rm T}\bm F_{i+1}$
		\EndFor
		\State let $\bm Z:=\bm I+\bm M_1^{\rm T}\bm F_{1,2:end}$\Comment{{\scriptsize $\bm F_{1,2:end}$ represents the second to last columns of $\bm F_1$.}}
		\State $\bm Z:=\bm Z^{-1}\bm M_1^{\rm T}$
		\State $\tilde{\bm \lambda}:=\bm F_{,1}-\bm F_{,2:end}\bm Z\bm F_{1,1}$\Comment{{\scriptsize $\bm F_{,1}$ represents the first column of $\bm F$, and $F_{1,1}$ represents the first column of $\bm F_1$.}}
	\end{algorithmic}
\end{algorithm}
\begin{remark}
	In (\ref{eqn:Sherman-Morrison-Woodbury formula}), the main cost is solving $(\bm L\bm \Lambda\bm L^{\rm T})^{-1}\tilde{\bm F}$ and $(\bm L\bm \Lambda\bm L^{\rm T})^{-1}\bm E_1$, which can be executed in parallel. In summary, the time complexity for solving (\ref{eqn:lagrangeSystemBlock}) is $O(Nn_m^2+(n_D+n_N)n_m^3)$, which is more effective than the directly $LDL^{\rm T}$ factorization when $N$ is large. Furthermore, the space complexity only depends on the number of different $\phi_{0,m}$ and is independent of $N$.
\end{remark}
\begin{remark}
	When $N$ is small, algorithm \ref{alg:Double-Newton} will cost most of the time during the whole FETI algorithm, and the required storage to solve (\ref{eqn:lagrangeMultiplierSystem}) directly is small. Therefore, in this situation, the ``backslash" operator in Matlab will be a better choice.
\end{remark}
\section{Numerical results}\label{sect: Numerical results}
In this section, we show the efficiency of the proposed algorithm by some numerical examples. Before introducing the numerical examples, let's give some default settings.
\begin{itemize}
	\item Material parameters:\\
	piezoelectric substrate: YX128LN;\\
	electrode: Al;
	\item Size of each subregion:\\
	piezoelectric substrate of each unit block: $1\mu m\times0.1\mu m\times10\mu m$;\\
	electrode: $0.5\mu m\times 0.1\mu m\times 0.15\mu m$;\\
	thickness of PML: $2\mu m$;
	\item The number of grid points in each coordinate direction of each subregion:\\
	piezoelectric substrate of each unit block: $17\times 2\times 17$;\\
	electrode: $9\times 2\times 5$;\\
	PML's thickness direction: $5$;
	\item Dimensionless parameters:\\
	$c_1=10^{10}$, $\omega_1=10^7$, $\varepsilon_1=10^{-10}$, $e_1=1$, $l_1=10^{-2}$ and $\rho_1=1$;
	\item The damping functions of each PML regions:\\
	$\Omega_{L}^{\rm PML}$: $d_1^{(L)}(x_1)=(1-(\frac{x_1-x^{L}_1}{d_{\rm PML}})^2)^2$, $d_2^{(L)}(x_2)=0$, $d_3^{(L)}(x_3)=0$;\\
	$\Omega_R^{\rm PML}$: $d_1^{(R)}(x_1)=(1-(\frac{x_1-x^{R}_1}{d_{\rm PML}})^2)^2$, $d_2^{(R)}(x_2)=0$, $d_3^{(R)}(x_3)=0$;\\
	$\Omega_B^{\rm PML}$: $d_1^{(B)}(x_1)=0$, $d_2^{(B)}(x_2)=0$, $d_3^{(B)}(x_3)=(1-(\frac{x_3-x^{B}_3}{d_{\rm PML}})^2)^2$;\\
	$\Omega_{LB}^{\rm PML}$: $d_1^{(LB)}(x_1)=(1-(\frac{x_1-x^{LB}_1}{d_{\rm PML}})^2)^2$, $d_2^{(LB)}(x_2)=0$, $d_3^{(LB)}(x_3)=(1-(\frac{x_3-x^{LB}_3}{d_{\rm PML}})^2)^2$;\\
	$\Omega_{RB}^{\rm PML}$: $d_1^{(RB)}(x_1)=(1-(\frac{x_1-x^{RB}_1}{d_{\rm PML}})^2)^2$, $d_2^{(RB)}(x_2)=0$, $d_3^{(RB)}(x_3)=(1-(\frac{x_3-x^{RB}_3}{d_{\rm PML}})^2)^2$,\\
	where $d_{\rm PML}$ is the thickness of each PML region and $x_i^m$ is the coordinates of the $i$-direction at the junction of $\Omega_m^{\rm PML}$ and $\Omega$.
\end{itemize}
All the numerical examples are discretized by the quadratic Lagrange FEM, and implemented in Matlab R2023a in an AMD Ryzen 7, 2.90GHz CPU (8 cores, 16 threads) on a laptop computer with a 16 GB RAM.
\subsection{Comparing the FETI algorithm with the quadratic Lagrange FEM}\label{sect:NumericalResult 1}
\begin{figure}[!htp]
	\centering
	\begin{subfigure}{\textwidth} 
		\centering
		\includegraphics[width=0.9\textwidth]{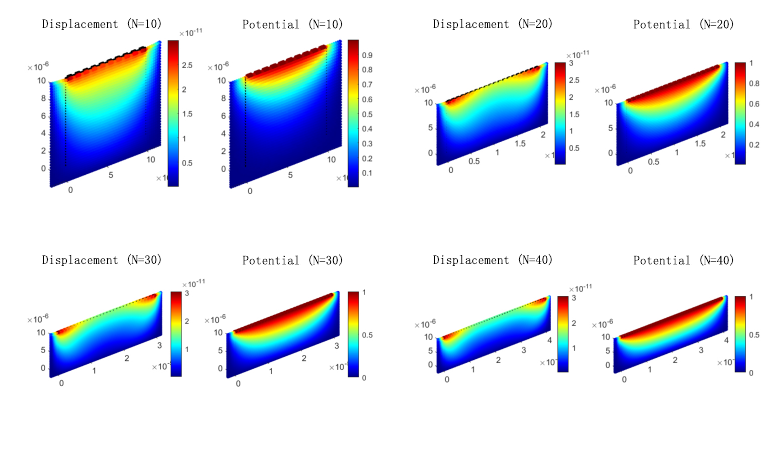} 
		\caption{}
	\end{subfigure}
	
	\begin{subfigure}{\textwidth} 
		\centering
		\includegraphics[width=0.9\textwidth]{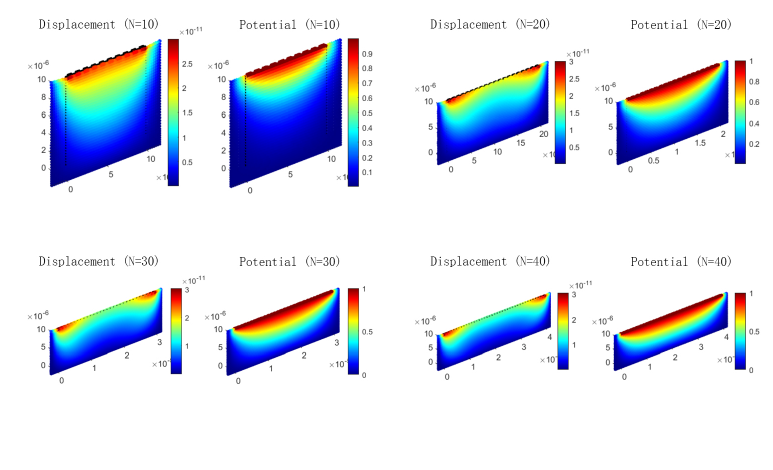} 
		\caption{}
	\end{subfigure}
	\caption{(a) Numerical solution of the FETI algorithm; (b) Numerical solution of FEM calculation (all deformations in the figure are magnified by $50000$ times).}
	\label{fig:FETIEX1Displacement}
\end{figure}

In the first numerical example, we consider the SAW equations with a voltage of $1V$ applied to each electrode, and solve it by the proposed method and the quadratic Lagrange FEM, respectively. We solve the SAW equations for $N=10, 20, 30, 40, 50$ using both methods. The runtime of the FETI and the FEM is shown in Tab \ref{Tab:EX1RunTime}, and Fig \ref{fig:FETIEX1Displacement}show the numerical solutions obtained by these two methods.
\begin{table}[H]
	\begin{tabular}{|l|llll|l|}
		\hline
		& \multicolumn{4}{c|}{FETI}                                                                                           & \multicolumn{1}{c|}{FEM}          \\ \hline
		$N$ & \multicolumn{1}{l|}{assemble  $\bm A$ and $\bm b$} & \multicolumn{1}{l|}{solve  $\bm A\bm \lambda=\bm b$} & \multicolumn{1}{l|}{compute $\bm X_m^p$, $\bm X_m^e$} & total runtime & total runtime \\ \hline
		10     & \multicolumn{1}{l|}{58.490}        & \multicolumn{1}{l|}{1.109}      & \multicolumn{1}{l|}{0.031}   & 59.630        & 25.861        \\ \hline
		20     & \multicolumn{1}{l|}{57.678}        & \multicolumn{1}{l|}{2.266}      & \multicolumn{1}{l|}{0.035}   & 59.979        & 53.329        \\ \hline
		30     & \multicolumn{1}{l|}{57.145}        & \multicolumn{1}{l|}{3.669}      & \multicolumn{1}{l|}{0.053}   & 60.867        & 120.031       \\ \hline
		40     & \multicolumn{1}{l|}{57.280}        & \multicolumn{1}{l|}{4.901}      & \multicolumn{1}{l|}{0.070}   & 62.251        & 250.989       \\ \hline
		50     & \multicolumn{1}{l|}{57.258}        & \multicolumn{1}{l|}{6.046}      & \multicolumn{1}{l|}{0.094}   & 63.398        & out of memory \\ \hline
	\end{tabular}
	\caption{the runtime (s) spends in the FEM and each step of the FETI algorithm.}
	\label{Tab:EX1RunTime}
\end{table}

\begin{figure}[H]
	\centering
	\includegraphics[width=10cm]{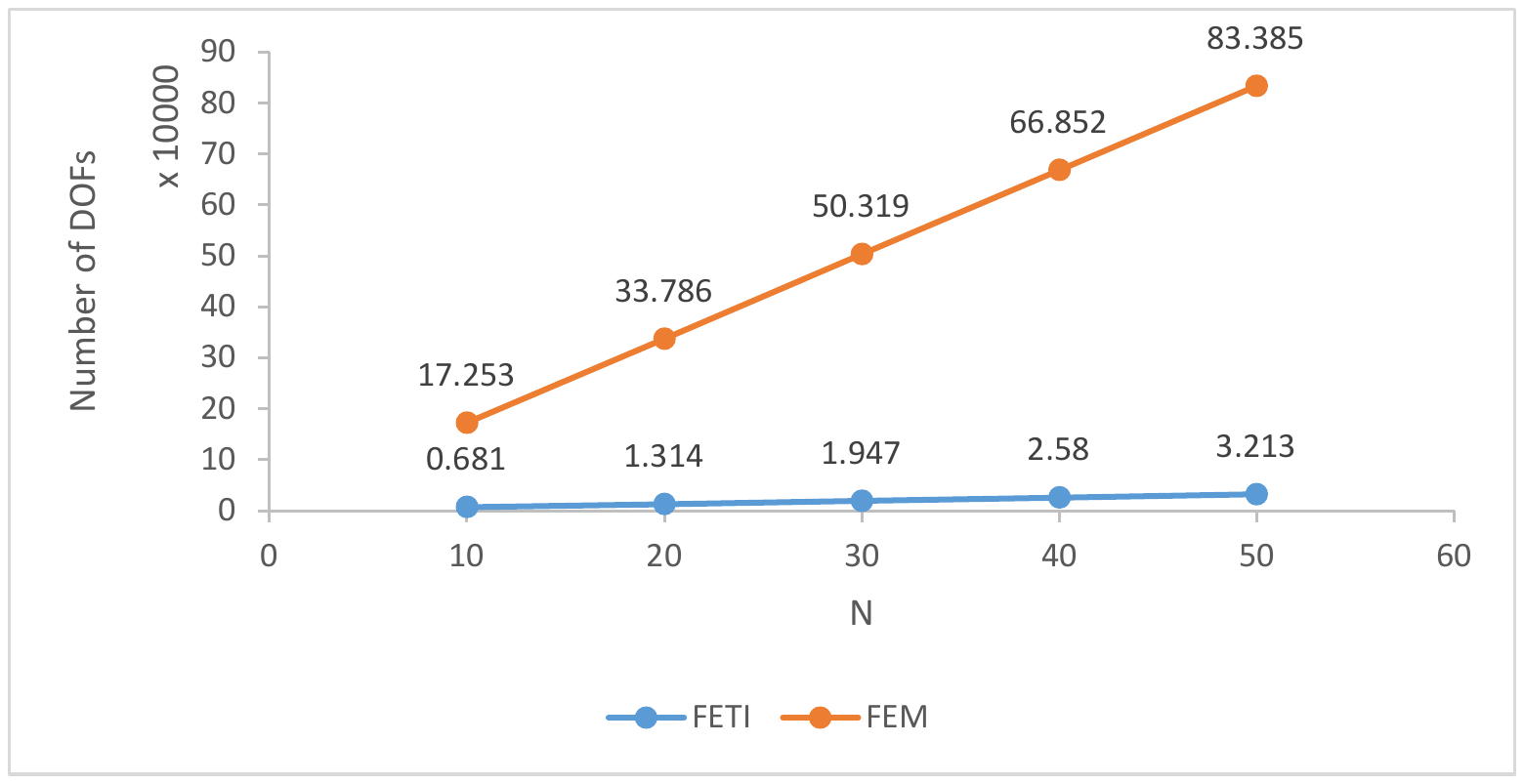}
	\caption{The DOFs of ``$\bm A\bm\lambda=\bm b$" and the FEM linear system}
	\label{fig:DOF}
\end{figure}

In Tab. \ref{Tab:EX1RunTime}, the linear system $\bm A\bm\lambda=\bm b$ is solved by the ``backslash" operator in Matlab, since the size of $\bm A$ is small. The process of ``assemble $\bm A$ and $\bm b$" is (\ref{eqn:matrixBlock})-(\ref{eqn:vectorBlock}), which is independent of $N$, and it acts as the main computational quantity of the FETI algorithm when $N$ is not very large.
The step ``computing $\bm X_m^p$, $\bm X_m^e$" (i.e. (\ref{eqn:numerical solution})) takes very little time, and it will be independent with $N$, if the number of threads of your computer is bigger than $N$. Moreover,  from the total runtime of FETI and FEM, we can observe that the advantage of the proposed algorithm becomes more and more obvious as $N$ increases. In this numerical example, only the size of $\bm A\bm\lambda=\bm b$ increases with $N$ during the entire FETI algorithm. Fig. \ref{fig:DOF} shows the size of this linear system and the sizes of the FEM's linear system as $N$ increases. \\
Fig. \ref{fig:DOF} shows that the required storage of the FETI algorithm is much smaller than that of the FEM, although we solve ``$\bm A\bm \lambda=\bm b$" directly by the ``backslash" operator in Matlab. This is also the reason why FEM runs out of memory when $N=50$, while the FETI algorithm can still obtain the numerical solution. In fact, the number of DOFs of  ``$\bm A\bm\lambda=\bm b$" is equal to that of all the subregions' interfaces.Therefore, the FETI algorithm has a dimensionality reduction effect when solving the SAW equations.
\subsection{Applying different voltages on different electrodes}\label{sect:NumericalResult 2}
\begin{figure}[H]
	\centering
	\includegraphics[width=15cm]{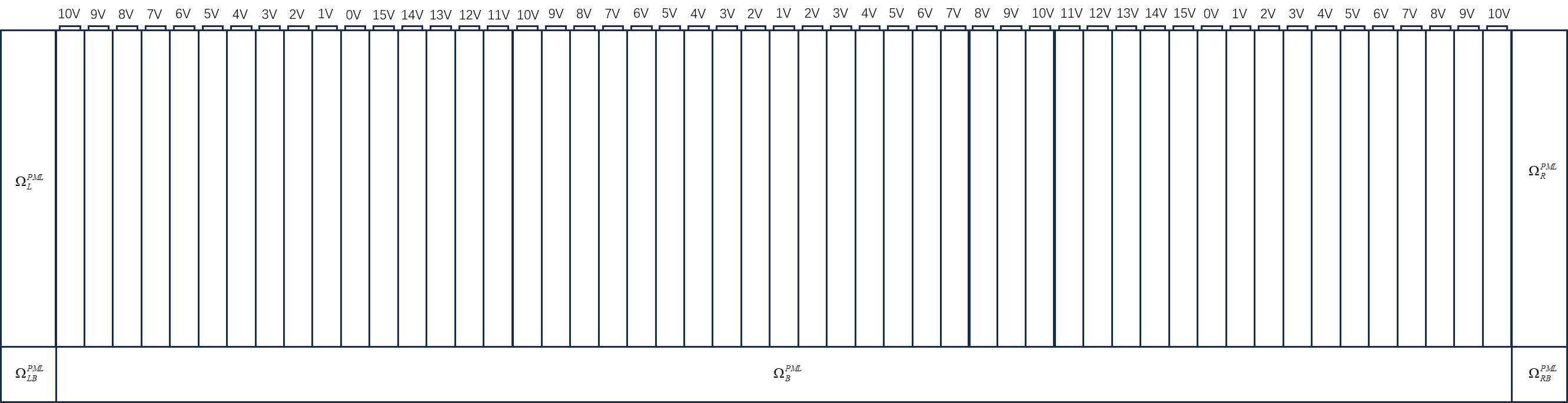}
	\caption{the voltages applied to each electrode with $j=15$.}
	\label{fig:51IDTs}
\end{figure}

Now, we fix $N=51$, and apply different voltages on the electrodes. The voltages $\phi_{0,i}^j$ are defined as follows.
$$
\phi_{0,i}^j=|i-25|\bmod j,\quad i=1,...,51,
$$
where $j$ is the number of different voltages. As an example, For example, the voltages applied to each electrode when $j=15$ are shown in the Fig. \ref{fig:51IDTs}.

\begin{figure}[H]
	\centering
	\includegraphics[width=10cm]{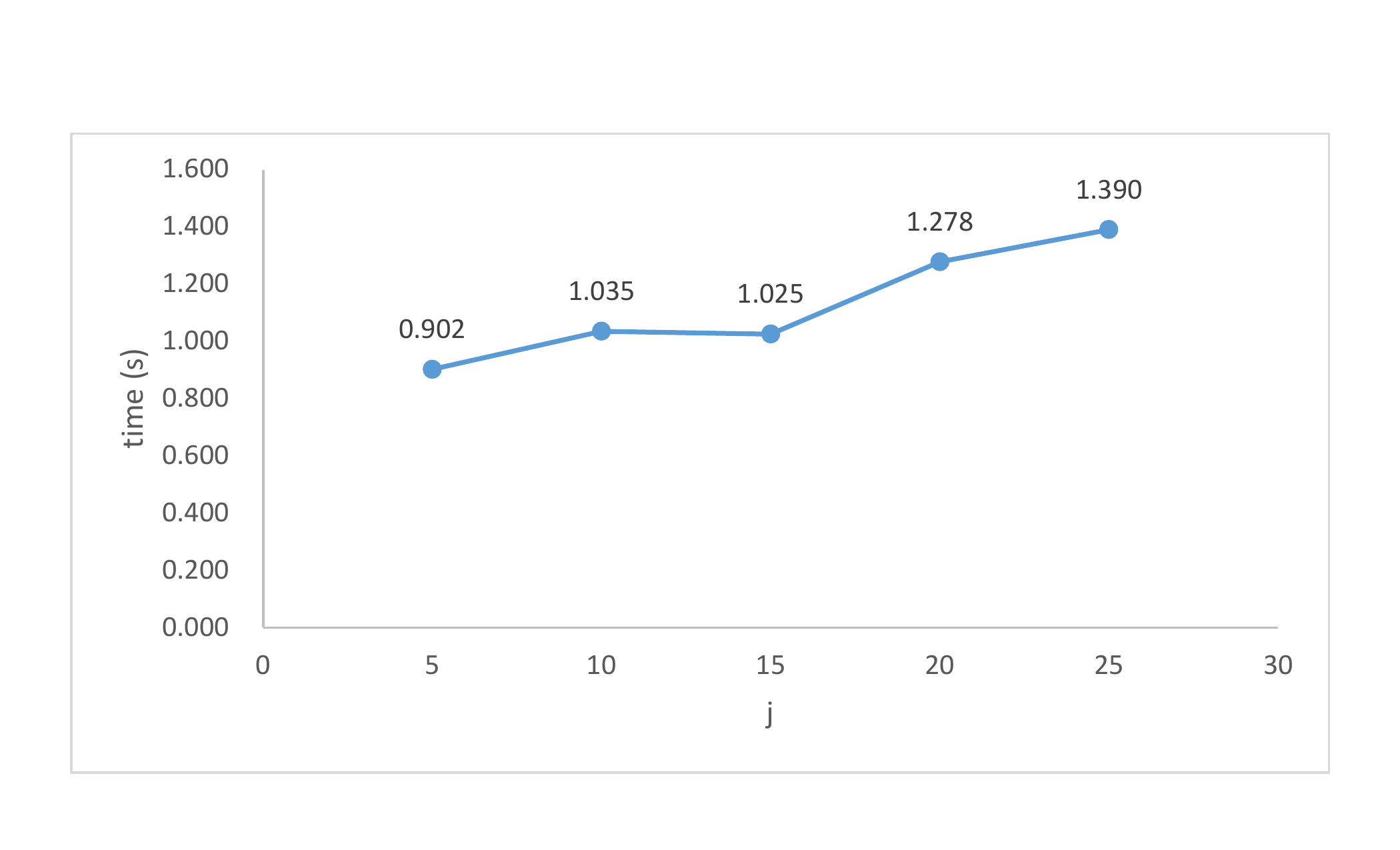}
	\caption{the required time of assembling $\bm b$ with different $j$.}
	\label{fig:TimeDifferentPotentials}
\end{figure}

\begin{figure}[!htp]
	\centering
	\includegraphics[width=1\textwidth]{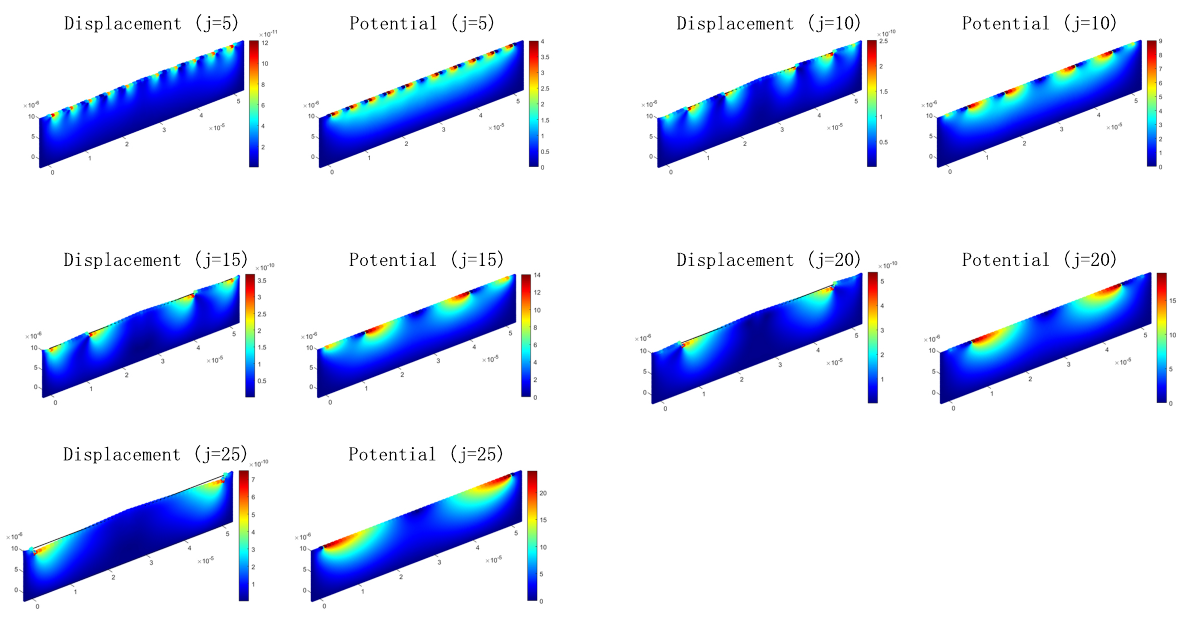}
	\caption{the numerical solutions of the SAW equations with different $j$ (The distortion in each figure has been magnified by $2000$ times).}
	\label{fig:DifferentPotentials}
\end{figure}

According to Sect. \ref{sect:FETI}, only step (\ref{eqn:vectorBlock}) is affected by $j$, since different voltages only cause variations in $\bm F_m^p$. Therefore, different values of $j$ only affect the time required to assemble the vector $\bm b$. Fig. \ref{fig:TimeDifferentPotentials} shows the time taken to assemble $\bm b$ for different $j$, and Fig. \ref{fig:DifferentPotentials} shows the corresponding numerical solutions.
In (\ref{eqn:vectorBlock}), since we can compute $\bm K_p^{-1}\bm F_m^p$  for the unit blocks in parallel with different $\phi_{0,i}^j$, applying different voltages on the electrodes will not significantly increase the runtime of the proposed FETI algorithm, if your computer has enough threads.

\subsection{Large-scale SAW equation}\label{sect:NumericalResult 3}
In the final numerical example, we consider the large-scale SAW equations. In this case, our laptop does not have enough memory to store the matrix $\bm A$, even if it could be stored, solving $\bm A\bm \lambda=\bm b$ by the ``backslash" in Matlab would take quite a bit of time. Therefore, we use the method proposed in Sect. \ref{Sect:Solve qusi-Toeplitz} to solve $\bm A\bm \lambda=\bm b$, during the execution of the FETI algorithm. First, let us show the efficiency and the accuracy of the proposed Double-Newton method.

We use the Double-Newton method, the Newton's method of (\ref{eqn:quadratic matrix equation}), the Double method and the Newton's method of (\ref{eqn:main nonlinear system}) to solve (\ref{eqn:main nonlinear system}) with initial value $\bm M$, respectively, and the iteration is terminated when the relative difference in the Frobenius-norm between the current step's solution and the previous step's solution is lower than $10^{-10}$. The results are shown in Tab. \ref{tab:Compare the 4 iter method}.
\begin{table}[H]
	\centering
	{\small\begin{tabular}{|l|l|l|l|l|}
		\hline
		& Double-Newton & Double   & Newton's method of (\ref{eqn:quadratic matrix equation}) & Newton's method of (\ref{eqn:main nonlinear system}) \\ \hline
		time (s)      & 68.765         & 1.327     & 77.543          & 64.035          \\ \hline
		$Err$       & 1.14E-11      & 3.50E-03 & 4.81E-02       & 3.27E-08       \\ \hline
		iterations & 16            & 8        & 8              & 7              \\ \hline
	\end{tabular}}
	\caption{\small Comparison of the four iterative methods to solve (\ref{eqn:main nonlinear system}).}
	\label{tab:Compare the 4 iter method}
\end{table}
In Tab. \ref{tab:Compare the 4 iter method}, ``iterations" refers to the number of iterations when the
iteration method is terminated, and the relative error $Err$ is defined as
$$
Err:=\frac{\|\bm B\bm \Lambda_{1*}^{-1}\bm B^{\rm T}+\bm \Lambda_{1*}-\bm M\|_F}{\|\bm M\|_F}.
$$
In Tab. \ref{tab:Compare the 4 iter method},
we can see that the double method terminates the fastest, but the accuracy is poor. The proposed Double-Newton method
 has the best accuracy.

\begin{figure}[!htp]
	\centering
	\begin{subfigure}{\textwidth} 
		\centering
		\includegraphics[width=1\textwidth]{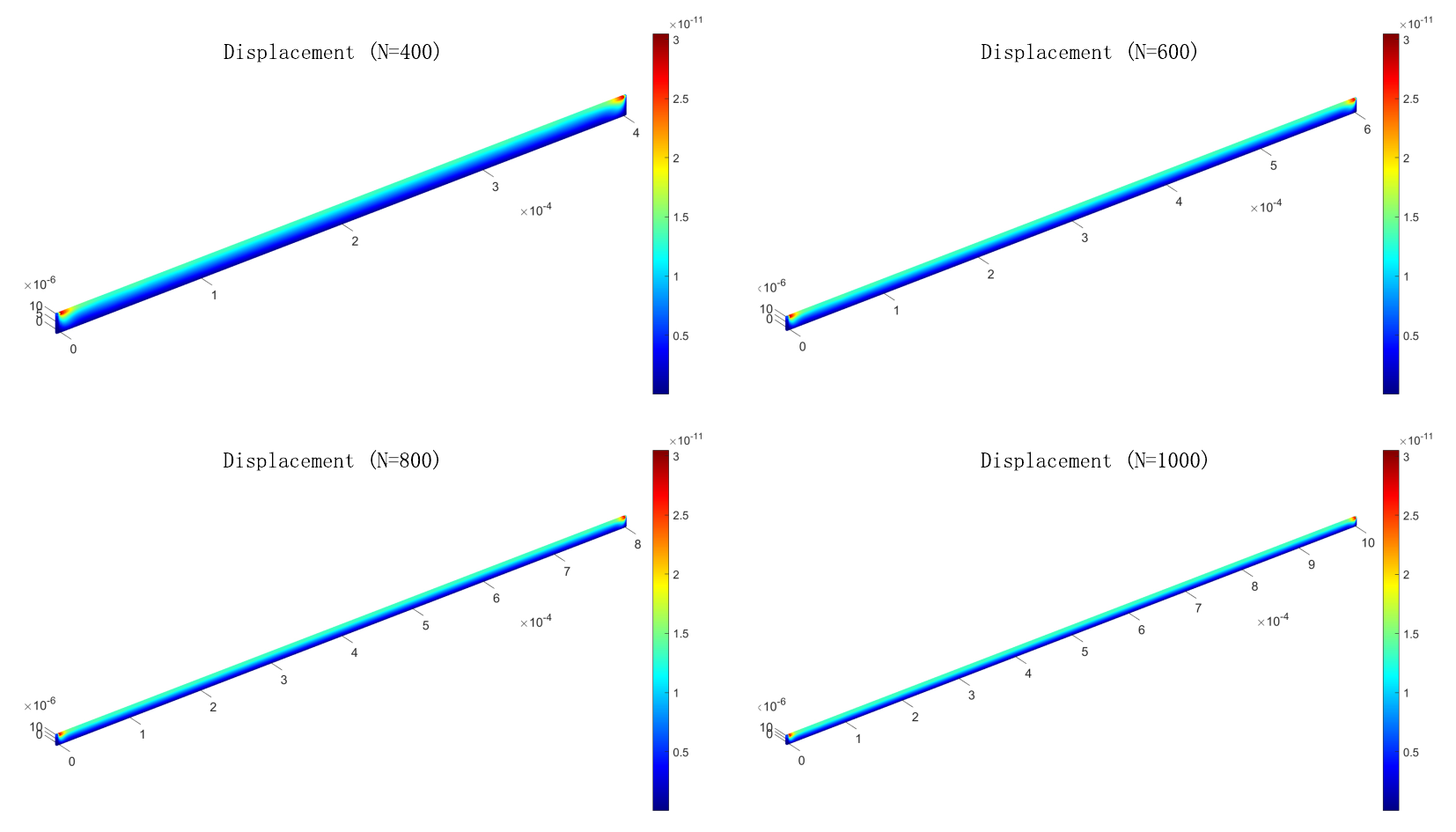} 
		\caption{}
	\end{subfigure}
	
	\begin{subfigure}{\textwidth} 
		\centering
		\includegraphics[width=1\textwidth]{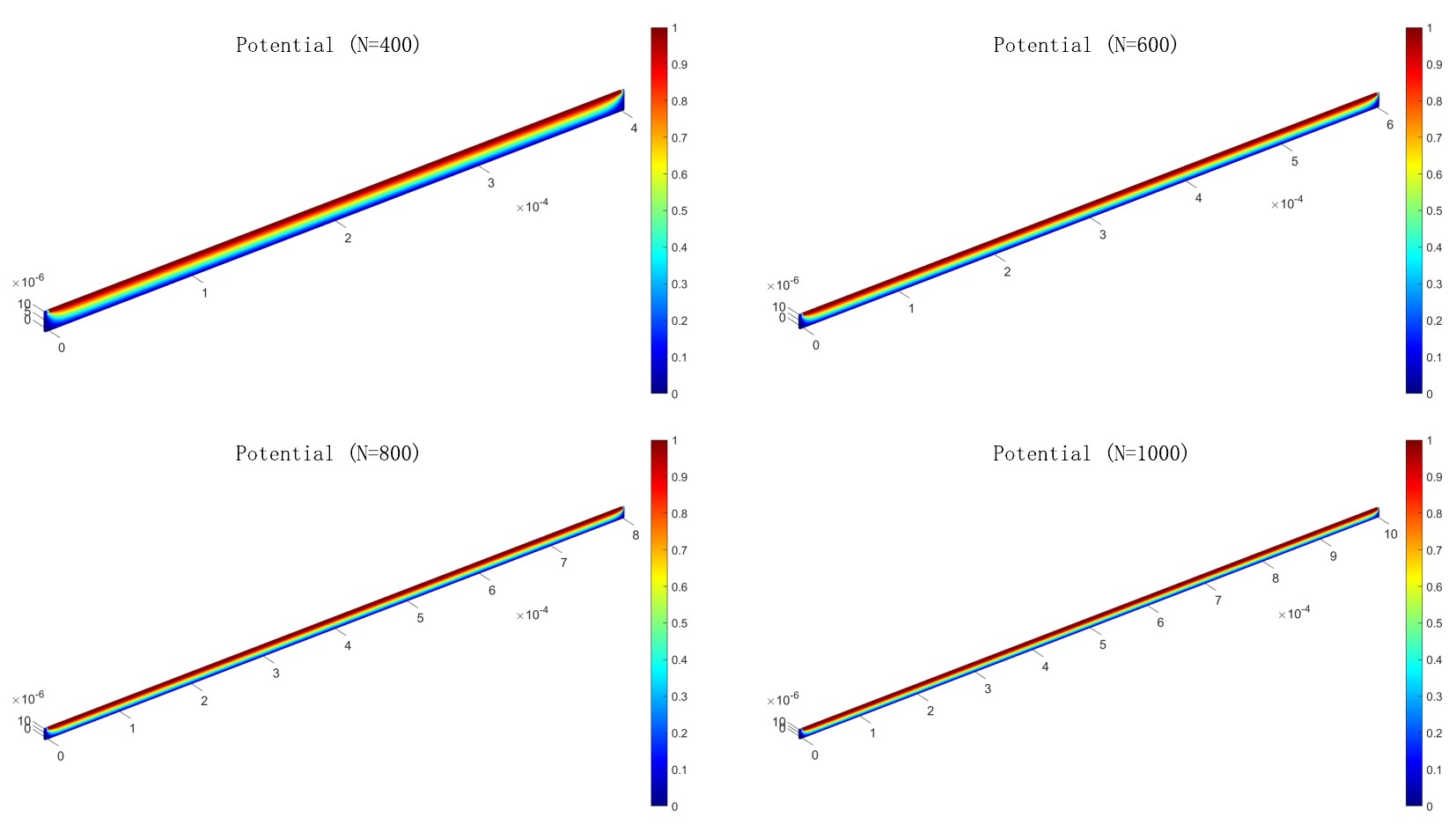} 
		\caption{}
	\end{subfigure}
	\caption{Numerical solutions of large-scale SAW equations: (a) Displacement (all deformations in the figure are magnified by 10000 times); (b) Electric potential.}
	\label{fig:FETIthe deformations of the large SAW device}
\end{figure}

Next, we use the FETI algorithm equipped with the Double-Newton method to solve the SAW equations for $N=400$, $600$, $800$ and $1000$. The time spent on each step is shown in Tab. \ref{Tab:The time (s) spent on each step}, and the numerical solutions of each SAW equations are shown in Fig. \ref{fig:FETIthe deformations of the large SAW device}.
\begin{table}[H]
	\centering
	{\begin{tabular}{|l|l|l|l|}
		\hline
		N    & assemble the matrix blocks & solve (\ref{eqn:matrixEquations})&compute (\ref{eqn:Sherman-Morrison-Woodbury formula})   \\ \hline
		400  & 41.324 & 68.588           & 30.537                          \\ \hline
		600  & 41.527 & 68.327           & 93.535                         \\ \hline
		800  & 40.976 & 68.539           & 223.420                         \\ \hline
		1000 & 40.686 & 68.189           & 246.373                         \\ \hline
	\end{tabular}}
	\caption{\small The time (s) spends on each step of solving $\bm A\bm \lambda=\bm b$.}
	\label{Tab:The time (s) spent on each step}
\end{table}
In Tab. \ref{Tab:The time (s) spent on each step}, the step ``assemble the matrix blocks" is the time spent on computing and assembling the matrix blocks  $\bm M$, $\bm M_L$, $\bm M_R$, $\bm B$, $\bm b_L$ and $\bm b_i$ defined in (\ref{eqn:MatrixBlocks}). From the table, we can see that only the time spent on computing (\ref{eqn:Sherman-Morrison-Woodbury formula}) increases with $N$. Tab. \ref{Tab:The time (s) spent on each step} also shows that the time required to solve (\ref{eqn:Sherman-Morrison-Woodbury formula}) grows slower and slower, when $N>600$. The reasons for this phenomenon is that the number of columns of $\hat{\bm F}_1$ exceeds that of $\bm \Lambda_{1}$ in line 6 of Algorithm \ref{alg:Sherman-Morrison-Woodbury formula}, when $N=800$ and $1000$. Therefore, we can compute $\bm \Lambda_{1}^{-1}\hat{\bm F}_1$ by first computing $\bm\Lambda_{1}^{-1}$ and then computing $\bm \Lambda_{1}^{-1}\hat{\bm F}_1$ instead of solving $n_c$ linear systems, where $n_c$ is the number of columns of $\hat{\bm F}_1$. As a result, the time complexity changes from $O(Nn_m^3)$ to $O(Nn_m^2)$.

When $N=400$, $600$, $800$ and $1000$, the total number of DOFs for the displacement and potential in the SAW equations are about $7.08$ million, $10.61$ million, $14.15$ million and $17.68$ million, respectively. It is almost impossible to solve them using the FEMs on a daily-used laptop. However, the FETI algorithm with the Double-Newton method proposed in this paper makes it possible to solve these problems on a laptop with only 16 GB of memory, and the computation time is also acceptable. This also reflects that the spatial complexity of the proposed algorithm is very low.
\section{Conclusion}
In this paper, we propose a FETI algorithm for solving the SAW equations defined in domains with periodic structures. Due to the periodic structure, the domain can be divided into four different types of subregions. This is the main advantage of solving this model using the FETI algorithm. In the first numerical example, the FETI algorithm performs more efficiently than the classical quadratic Lagrange FEM when $N$ becomes large. We also consider the SAW equations with different voltages applied to the electrodes, and the numerical results show that it does not significantly reduce the efficiency of the FETI algorithm.

The related large scale linear system corresponding to the Lagrange multipliers is efficiently solved by a direct solver. The computational complexity is very low and the computation time is also reasonable. This approach overcomes the storage challenges posed by the large problem size. The numerical results show that the proposed method performs efficiently

For the sake of presentation, this paper only considers the case where the right-hand load vanishes. However, the proposed method can be extended to the case with periodic load, as shown in the second numerical example.




{
	\rm
	\bibliographystyle{plain}
	\bibliography{Reference}
}

\end{document}